\UseRawInputEncoding

\documentclass[12pt, reqno]{amsart}
\usepackage{amssymb,amsthm,amsmath}
\usepackage[numbers,sort&compress]{natbib}
\usepackage{amssymb,amsmath}
\usepackage{amsfonts}
\usepackage{mathrsfs}
\usepackage{latexsym}
\usepackage{amssymb}
\usepackage{amsthm}
\usepackage{pdfsync}
\usepackage{indentfirst}
\usepackage{xcolor}
\date{June 18, 2021}

\usepackage{amsmath}
\usepackage{amsthm}

\newtheorem{theorem}{Theorem}

\newtheorem{remark}[theorem]{Remark}
\newtheorem{proposition}[theorem]{Proposition}
\newtheorem{lemma}[theorem]{Lemma}
\newtheorem{definition}[theorem]{Definition}
\newtheorem{corollary}[theorem]{Corollary}

\newcommand{\beq}{\begin{equation}}
\newcommand{\eeq}{\end{equation}}
\newcommand{\ben}{\begin{eqnarray}}
\newcommand{\een}{\end{eqnarray}}
\newcommand{\beno}{\begin{eqnarray*}}
\newcommand{\eeno}{\end{eqnarray*}}

\numberwithin{equation}{section}

\allowdisplaybreaks

\begin{document}
\title[Partial regularity of the 3D chemotaxis-Navier-Stokes equations]{\bf   Partial regularity of solutions to the 3D chemotaxis-Navier-Stokes equations at the first blow-up time}
\author{Xiaomeng~Chen}
\address[Xiaomeng~Chen]{School of Mathematical Sciences, Dalian University of Technology, Dalian, 116024,  China}
\email{cxm2381033@163.com}

\author{Shuai~Li}
\address[Shuai~Li]{School of Mathematical Sciences, Dalian University of Technology, Dalian, 116024,  China}
\email{leeshy@mail.dlut.edu.cn}

\author{Wendong~Wang}
\address[Wendong~Wang]{School of Mathematical Sciences, Dalian University of Technology, Dalian, 116024,  China}
\email{wendong@dlut.edu.cn}
\date{\today}
\maketitle

\begin{abstract}As Dombrowski et al. showed in \cite{DCCGK2004} (see also \cite{TCDWKG2005}), suspensions of aerobic bacteria often develop flows from the interplay of chemotaxis and buoyancy, which is described as the chemotaxis-Navier-Stokes flow, and they observed self-concentration occurs as a turbulence by exhibiting transient, reconstituting, high-speed jets. Moreover, local concentration leads to a jet descending faster than its surroundings, which entrains nearby fluid to produce paired, oppositely signed vortices.
In this note, we investigate the  Hausdorff dimension of these vortices (singular points) by considering partial regularity of weak solutions of the three dimensional  chemotaxis-Navier-Stokes equations, and obtain the $\frac53$-dimensional Hausdorff measure of the possible singular set is vanishing  at the first blow-up time, which generalizes the Caffarelli-Kohn-Nirenberg's partial regularity theory to the chemotaxis-fluid model. The new ingredients  are to establish certain type of local energy inequality and deal with the non-scaling invariant quantity of $n\ln n$, where $n$ represents the cell concentration, which seems to be the first description for the singular set of weak solutions  of the model.
\end{abstract}

{\small {\bf Keywords:} chemotaxis-Navier-Stokes, partial regularity, local energy inequality}

{\bf 2010 Mathematics Subject Classification:} 35Q30; 35Q35; 76D05.

\setcounter{equation}{0}
\section{The question}

Consider a PDE model on $Q_T=\mathbb{R}^3\times(0,T)$ describing the dynamics of oxygen,
swimming bacteria, and viscous incompressible fluids, which was proposed by Tuval et al. \cite{TCDWKG2005} as follows:
\begin{eqnarray}\label{eq:GKS}
 \left\{
    \begin{array}{llll}
    \displaystyle \partial_t n+u\cdot \nabla n-\Delta n=-\nabla\cdot(\chi(c)n\nabla c),\\
    \displaystyle \partial_t c+ u\cdot \nabla c-  \Delta c=-\kappa(c)n, \\
    \displaystyle \partial_t u+ u\cdot \nabla u-\Delta u+\nabla p =-n\nabla \phi,~~\nabla\cdot u=0 \\
    \end{array}
 \right.
\end{eqnarray}
where $c(x,t):Q_T\rightarrow{\mathbb R}^{+}$, $n(x,t):Q_T\rightarrow\mathbb{R}^{+}$, $u(x,t):Q_T\rightarrow\mathbb{R}^{3}$ and $p(x,t):Q_T\rightarrow\mathbb{R}$ denote the oxygen concentration, cell concentration, the fluid velocity and the associated pressure, respectively. Moreover, the
gravitational potential $\phi$, the chemotactic sensitivity $\chi(c)\geq 0$ and the per-capita oxygen
consumption rate $\kappa(c)\geq 0$ are sufficiently smooth given functions (see also \cite{DCCGK2004}).

Due to the significance of the biological background (see \cite{DCCGK2004}, \cite{TCDWKG2005}), the model could be used to predict the
large-scale bioconvection affecting clearly the overall oxygen consumption in the
above experiments (see, for example, \cite{AL2010}). Many mathematicians have studied this model and made much progress, such as the existence of weak solutions, the chemotaxis-Navier-Stokes system with a nonlinear diffusion,  blow-up criteria, stability and so on. Here we just mention some related works for the result in this paper.

Firstly, for the existence of weak solutions,  in \cite{DLM2010}, global classical solutions near constant steady states are constructed for the full
chemotaxis-Navier-Stokes system by Duan-Lorz-Markowich. In \cite{AL2010}, for the case of bounded domain in $\mathbb{R}^n$ with $n= 2,3$, the local existence of weak solutions for problem (\ref{eq:GKS}) is obtained by Lorz. Later, Winkler proved the existence of global weak solution in \cite{Winkler2012} by assuming that
\beno
\left(\frac{\kappa}{\chi}\right)'>0,\quad \left(\frac{\kappa}{\chi}\right)''\leq 0,\quad  \left({\kappa}{\chi}\right)'\geq0.
\eeno
By assuming $\chi',\kappa'\geq 0$ and $\kappa(0)=0$. In \cite{CKL2013} and \cite{CKL2014}, local well-posed results and blow-up criteria were established by Chae-Kang-Lee. Recently, Winkler proved the global existence of weak solutions of the system (\ref{eq:GKS}) in bounded domain with large initial data, and obtianed much better a priori estimates such as $\frac{|\nabla c|^4}{c^3}\in L^1$ in \cite{Winkler2016}. For the two-dimensional system  of (\ref{eq:GKS}), the system is better understood. Liu and Lorz \cite{LL2011}  proved the global existence of weak solutions to the two-dimensional system  of (\ref{eq:GKS}) for arbitrarily large initial data, under the assumptions on $\chi$ and $f$ made in \cite{DLM2010}. See \cite{LL2016,DLX2017,WWX2018,WWX2018-2,He2020,WZZ2021} and the references therein for more results.
For more references  about the existence of solutions, we refer to \cite{CL2016,Winkler2017,KM2019,DL2020,Bl2020,DL2022} and the references therein.

As for the case of the chemotaxis-Navier-Stokes system with a nonlinear diffusion, that means $\Delta n$ is replaced by $\Delta n^m$, there are also many results.
In \cite{FLM2010}, Lorz-Francesco-Markowich showed the global existence of a bounded solution to porous medium equation, on a bounded domain in $\mathbb{R}^2$, with the boundary conditions $\partial_v n^m=\partial_v c=u=0$ and the condition $m\in (\frac32,2]$. In \cite{TW2012}, Tao and Winkler extended the result to the case $m > 1$ on a bounded domain in $\mathbb{R}^2$. In \cite{LL2011}, Lorz and Liu proved the global existence of a weak solution to porous medium equation in $\mathbb{R}^3$ when $m =\frac43$. For more references, one can refer to \cite{TW2013,CKK2014,ZK2021} and so on.

For the consideration of boundary conditions, inhomogeneous Dirichlet boundary conditions for the signal may affect global regularity in the three-dimensional full Navier-Stokes version, which is founded by  Black-Winkler in \cite{BW2022}, we referred the recent result.

As Winkler said in \cite{Winkler2016},
{ ``For the full three-dimensional chemotaxis-Navier-Stokes system, even at the very basic level of global existence in generalized solution frameworks, a satisfactory solution theory is entirely lacking.''} In this paper our aim is to explore partial regularity properties of weak solutions.
For simplicity, we consider the case $\kappa(c)=c$ and $\chi(c)=1$. Then the three dimensional chemotaxis-Navier-Stokes system (\ref{eq:GKS}) is reduced to
\begin{eqnarray}\label{eq:KS}
 \left\{
    \begin{array}{llll}
    \displaystyle \partial_t n+u\cdot \nabla n-\Delta n=-\nabla\cdot(n\nabla c),\\
    \displaystyle \partial_t c+ u\cdot \nabla c- \Delta c=-cn, \\
    \displaystyle \partial_t u+u\cdot \nabla u-\Delta u+\nabla p =-n\nabla \phi,~~\nabla\cdot u=0, \\
    \end{array}
 \right.
\end{eqnarray}
Recently, global weak solution for this system was obtained in 2D and 3D (see, for example, \cite{ZZ2014,HZ2017}), respectively, where they established a priori estimate
\beno
\mathcal{U}(t) + \int_0^t \mathcal{V}(t) d\tau \leq C e^{Ct},
\eeno
where
\beno
\mathcal{U} = \|n\|_{L^1 \cap L \log L} + \|\nabla \sqrt c\|_2^2 + \|u\|_2^2,
\eeno
and
\beno
\mathcal{V} = \|\nabla \sqrt{n+1}\|_2^2 + \|\Delta \sqrt c\|_2^2 + \|\nabla u\|_2^2 + \int_{\mathbb{R}^d} (\sqrt{c})^{-2} |\nabla \sqrt c|^4 dx + \int_{\mathbb{R}^d} n |\nabla \sqrt c|^2 dx,
\eeno
where $d=2,3$, the definition of $ L \log L$ is given by Definition  \ref{def:log}, and  we write the norm of $\|f\|_{L^q(\mathbb{R}^d)}$ as $\|f\|_q$ for simplicity. We also refer to the recent existence result of global weak solution in $\mathbb{R}^3$ in \cite{KLW2022} by assuming $(n_0+1)\ln(n_0+1)\in L^1$. They established a different priori estimate as following:
\beno
\frac d{dt} F_\varepsilon(t) + \frac12 D_\varepsilon(t) \leq C.
\eeno
Here,
\beno
F_\varepsilon(t) = \int_{B_{\frac1\varepsilon}} (n_\varepsilon+1) \ln (n_\varepsilon+1) + \frac12 \int_{B_{\frac1\varepsilon}} \frac{|\nabla c_\varepsilon|^2}{c_\varepsilon} + \frac b2 \int_{B_{\frac1\varepsilon}} |u_\varepsilon|^2,
\eeno
and
\beno
D_\varepsilon(t) &=& \frac14 \int_{B_{\frac1\varepsilon}} \frac{|\nabla n_\varepsilon|^2}{n_\varepsilon+1} + K_1 \int_{B_{\frac1\varepsilon}} \frac{|D^2 c_\varepsilon|^2}{c_\varepsilon} + \frac{K_1}4 \int_{B_{\frac1\varepsilon}} \frac{|\nabla c_\varepsilon|^4}{c_\varepsilon^3} \\
&&+ \frac12 \int_{B_{\frac1\varepsilon}} \frac{F_\varepsilon(n_\varepsilon)}{c_\varepsilon} |\nabla c_\varepsilon|^2 + \frac b2 \int_{B_{\frac1\varepsilon}} |\nabla u_\varepsilon|^2,
\eeno
with $F_\varepsilon(s) = \varepsilon^{-1} \ln (1+\varepsilon s)$.
However, up to now more information about these weak solutions is still not known,  especially the interior singular vortices as described in \cite{DCCGK2004} and the self-organized generation of a persistent hydrodynamic vortex that traps cells near the contact line(see \cite{TCDWKG2005}). Motivated by the recent progress on the non-uniqueness of suitable Leray-Hopf solutions to the Navier-Stokes equations with identical body force by Albritton-Bru\'{e}-Colombo in \cite{ABC}, it's interesting to consider suitable weak solutions of (\ref{eq:KS}) as Caffarelli-Kohn-Nirenberg in \cite{CKN}, and one may ask naturally:
\\
{\it Q1: Whether does there exist a suitable weak solution for the system of (\ref{eq:GKS}) or (\ref{eq:KS})? }\\
{\it Q2:  How to characterize the singular points of weak solutions?}

In this note we aim to answer the second question.
Recall that these so-called partial regularity or $\varepsilon$-regularity theory, it can be traced back to
the well-known work by Caffarelli-Kohn-Nirenberg \cite{CKN} for the analysis of suitable weak solutions of the three dimensional time-dependent Navier-Stokes equations, where they  showed that the set $\mathcal{S}$ of possible interior singular points of a suitable weak solution is one-dimensional parabolic Hausdorff measure zero by improving Scheffer's results in \cite{SV1,SV2,SV4}. The suitable weak solution is better than Leray-Hopf weak solution introduced by Leray in
\cite{Leray} and  if the local strong solution blows up, then the solution may be continued as a
suitable weak solution (see Proposition 30.1 in \cite{LR}).

Besides, it is worth mentioning the interesting approach
of Katz and Pavlovi\'{c} (\cite{KP2002}) for studying the dimension of the singular set, where they considered the Navier Stokes equation with dissipation $(-\Delta)^\alpha$
with the condition of $ 1<\alpha < \frac54$ and proved the Hausdorff dimension of the singular set at time of first blow up is at most $5-4\alpha$.
More references on simplified proofs and improvements, we refer to Lin \cite{Lin}, Ladyzhenskaya-Seregin \cite{LS}, Tian-Xin \cite{TX}, Seregin \cite{Se}, Gustafson-Kang-Tsai \cite{GKT}, Vasseur \cite{Va} and the references therein.
Here we consider the partial regularity of the system (\ref{eq:KS}) at the first blow-up time as Dong-Du in \cite{DD}, since the first question of $({\it Q1})$ is still unknown, which is an open question.

Recall the well-posed results in \cite{CKL2013} or \cite{HZ2017}.
\begin{theorem}\label{thm:localwellposed}
Assume that $n_0\geq 0, c_0\geq 0$ and $\nabla^k \phi\in L^\infty$ with $1\leq k\leq m$.
There exists a constant $T^\ast$, the maximal existence time, which depends on the norm of initial data, such that for any  $t<T^\ast$, if the initial data $(n_0,c_0,u_0)\in H^{m-1}{(\mathbb{R}^3)}\times H^m{(\mathbb{R}^3)} \times H^m{(\mathbb{R}^3)}$ with $m\geq 3$ satisfy  (\ref{eq:KS}) in $\mathbb{R}^3\times(0,T^\ast)$, then there exists a unique regular solution $(n, c, u)$ of (\ref{eq:KS}) satisfying $n\geq 0,c\geq 0$ and
\ben\label{ine:assumption}
(n, c, u)\in L^\infty (0, t;H^{m-1}{(\mathbb{R}^3)}\times H^m{(\mathbb{R}^3)}\times H^m{(\mathbb{R}^3)}),\nonumber\\\nonumber
(\nabla n, \nabla c, \nabla u)\in L^2 (0, t;H^{m-1}{(\mathbb{R}^3)}\times H^m{(\mathbb{R}^3)}\times H^m{(\mathbb{R}^3)}),\\\nonumber
(\partial_t n, \partial_tc,\partial_t u)\in L^\infty (0, t;H^{m-1}{(\mathbb{R}^3)}\times H^m{(\mathbb{R}^3)}\times H^m{(\mathbb{R}^3)}).
\een
\end{theorem}
Firstly, we give the definition of $L\log L$ norm, which will be used in the following time.
\begin{definition}\label{def:log}
The Zygmund classes with $A(t) = t\log^{+} t$, is defined as the set all functions $f$ such that
\beno
\int_{\mathbb{R}^3}A(|f(x)|)dx<\infty.
\eeno
the corresponding Zygmund space $L\log L(\mathbb{R}^3)$ is defined as the linear hull of the Zygmund class,
which is equipped with the Luxemburg norm
\beno
||f||_{L\log L}=\inf\left\{k \big|\int_{\mathbb{R}^3}A(\frac fk)dx\leq 1\right\}.
\eeno
and
\begin{eqnarray}\label{eq:initial}
 \log^{+}t=\left\{
    \begin{array}{llll}
    \displaystyle \log t,~~t\geq 1 ,\\
    \displaystyle 0,~~{\rm otherwise}.
    \end{array}
 \right.
\end{eqnarray}
\end{definition}

For the given initial data, there exists a global weak solution (see \cite{Winkler2016} or \cite{HZ2017}), which is defined as follows:
\begin{definition}\label{def:weaksolution}
$(n,c,u)$ is called a weak solution to the Cauchy problem (\ref{eq:KS}) with the initial data $(n_0,c_0,u_0)$ satisfying
\beno
n_0\in L^1\cap L \log L, n_0>0, c_0\in L^1\cap L^\infty, \sqrt{c_0}\in L^2, u_0\in L^2, {\rm div} u_0=0
\eeno
and $\nabla\phi\in L^\infty$, if the following conditions hold:\\

(i)  $n(t,x)>0$, $c(t,x)>0$,  for $t> 0$ and $x\in \mathbb{R}^3$, \\

(ii) $(n, c, u)$ satisfies the system (\ref{eq:KS}) in the sense of distribution;\\

(iii) $(n,c,u)$ satisfies: for any $t>0$, the following inequality is true
\beno
\mathcal{U}(t)+\int_0^t\mathcal{V}(\tau)d{\tau} \leq Ce^{Ct};
\eeno
where
\beno
&&\mathcal{U}(t)=\|n\|_{L^1\cap LlogL}+\|\nabla\sqrt{c}\|^2_{2}+\|u\|^2_2;\\
&&\mathcal{V}(t)=\|\nabla\sqrt{n+1}\|^2_2+\|\nabla^2\sqrt{c}\|^2_2+\|\nabla u\|^2_2+\int_{\mathbb{R}^3}(\sqrt{c})^{-2}|\nabla\sqrt{c}|^4dx+\int_{\mathbb{R}^3}n|\nabla\sqrt{c}|^2dx
\eeno
\end{definition}

%
%

Motivated by \cite{CKN}, we consider the partial regularity property of weak solutions, which is so-called $\varepsilon$-regularity criteria. First, we say
a point $(x_0,t_0)$ is a regular point if $(n,\nabla c,u)\in L^\infty(Q_{r_0}(x_0,t_0))$ for some $r_0>0$, where $(Q_{r_0}(x_0,t_0))=B_{r_0}(x_0)\times (t_0-r_0^2,t_0)$ and $B_{r_0}(x_0)=\{x, |x-x_0|<r_0\}$. It is worth noting that the definition here is consistent with the global regularity criterion proved in \cite{CKL2014}, where the regularity was ensured by
\beno
u\in L^2(0,T; L^\infty(\mathbb{R}^3)),\quad n\in L^1(0,T; L^\infty(\mathbb{R}^3)),
\eeno
(See Theorem 1 in \cite{CKL2014}). When $(x_0,t_0)=(0,0)$, we write $(0,0)=0$, $Q_{r_0}(0)=Q_{r_0}$ and $B_{r_0}(0)=B_{r_0}$ for simplicity.


Our first theorem is as follows:
\begin{theorem}\label{thm:regularity-fixed}
Assume that $(n,c,u)$ is a regular solution of (\ref{eq:KS}) in $\mathbb{R}^3\times(-1,0)$ with the initial data $(n(x,-1),c(x,-1),u(x,-1))\in  H^{2}{(\mathbb{R}^3)}\times H^3{(\mathbb{R}^3)} \times H^3{(\mathbb{R}^3)}$ as in Theorem \ref{thm:localwellposed}, which is also a weak solution as in Definition \ref{def:weaksolution}.  Then $z_0=(x_0,0)$ is a regular point, if there exists an absolute constant $\varepsilon_1$ such that
\ben\label{eq:condition-fixed}
&&\sup_{-1<t<0}\int_{B_1(x_0)} n + |n \ln n| + |\nabla \sqrt{c+1}|^2 + |u|^2dx \nonumber\\
&&+ \int_{Q_1(z_0)} |\nabla \sqrt{n}|^2 + |\nabla u|^2 + |\nabla^2 \sqrt{c+1}|^2 + |p|^\frac32dxdt \leq \frac{\varepsilon_1}{(\Lambda_0\Lambda_1)^{4+4\alpha_0}},
\een
where $\alpha_0>0$ is an absolute constant, $\Lambda_0=108\|c(\cdot,-1)+1\|_{L^\infty(\mathbb{R}^3)}$ and $\Lambda_1=(\|\nabla\phi\|_\infty+1).$
\end{theorem}

\begin{remark}\label{rem:1}
(i) The term of $\nabla \sqrt{n}\in L^2$ is reasonable, see the a priori estimates of weak solutions by Winkler in \cite{Winkler2016}, or  it can also  be derived from Lemma \ref{Local energy inequality} and the definition of weak solutions. The pressure $p$ is well-defined due to the Calder\'{o}n-Zygmund estimates and the equations of (\ref{eq:KS}). Especially, there holds
\beno
\int_{-1}^0\int_{\mathbb{R}^3}|p|^{\frac53}dxdt&\leq& C\int_{-1}^0\int_{\mathbb{R}^3}|u|^{\frac{10}3}dxdt+C\left(\int_{-1}^0\left(\int_{\mathbb{R}^3}|n\nabla\phi|^{\frac{15}{14}}dx\right)^{\frac{14}{9}}dt\right)\\
&\leq& C\int_{-1}^0\int_{\mathbb{R}^3}|u|^{\frac{10}3}dxdt+C\Lambda_1\|\sqrt{n}\|_{L^\infty_tL^2_x}^{3}\|\sqrt{n}\|_{L^2_tL^6_x}^{\frac13}\\
\eeno
since $u,\sqrt{n}\in L^\infty_tL_x^2$ and $\nabla u,\nabla\sqrt{n}\in L^2_tL_x^2$ imply that
\beno
u,\sqrt{n}\in L^s_tL_x^q, \quad \frac{2}{s}+\frac{3}{q}=\frac32, \quad 2\leq q\leq 6.
\eeno

(ii) The new observation of this theorem is the local energy inequality of (\ref{eq:local energy inequality}) (See Lemma \ref{Local energy inequality}), which is indeed a local a priori estimate for weak solutions, which is of independent interest.

(iii) The difficulty mainly lies in dealing with the term including $\ln n$, which is not scaling invariant under the embedding inequality.
 We establish the local a priori estimates of $\int_{B_1} (n \ln n \psi)(\cdot,t)$ firstly by estimating the local energy inequality, then use the equation of $n$ by
  estimating the term of $\int_{B_1} (n \psi)(\cdot,t)$, which combined with the embedding inequality in a fixed sphere imply the estimate of $\int_{B_1} n |\ln n| $.

(iv) In \cite{CKL2013,CKL2014}, the estimate of the term $\int n |\ln n|$ in whole space is that
\ben\label{whole}
\int n |\ln n| \leq \int n \ln n + 2 \int n (\ln n)_{-} \leq \int n \ln n + C + C \int n \langle x \rangle.
\een
Here $(\ln n)_{-}$ is a negative part of $\ln n$ and $\langle x \rangle^2 = 1 + |x|^2$. We use a different method to deal with this term without the weight $\langle x \rangle$, since the following estimate holds locally:
\ben\label{boundedness 1}
\int_{B_{r_k}} n |\ln n| \leq \int_{B_{r_k}} n \ln n + 2 \int_{B_{r_k}} n (\ln n)_{-} \leq \int_{B_{r_k}} n \ln n + C \int_{B_{r_k}} n^{1-\alpha}.
\een
\end{remark}

\begin{remark}To the authors' best
knowledge, whether
the weak solution considered as in \cite{HZ2017} verifies the local energy estimate is still unknown, and it seems that
the existence of such weak solutions  is still an
open problem. The main obstacle lies in the right hand term of $-\nabla\cdot(n\nabla c)$, since  the nonlinear term
could not be  controlled by the energy norm of $n$ and $\nabla\sqrt{c+1}$ under the Sobolev imbedding theorem in the energy inequality. Here we prove a local energy inequality of weak type with the term of $n\ln n$, which may have an uncertain sign.
\end{remark}

Similar as Lin's regularity version in \cite{Lin}, we have the following interior regularity criteria.
\begin{theorem}\label{thm:lin}
Taking the same assumptions as Theorem \ref{thm:regularity-fixed}, there exists a constant $\varepsilon_2(\Lambda_0,\Lambda_1) = \frac{\varepsilon_1^3}{C(\Lambda_0\Lambda_1)^{15+12\alpha_0}}$ and $\varepsilon_2'(\Lambda_0, \Lambda_1) = \frac{\varepsilon_1^\frac{15}4}{C(\Lambda_0\Lambda_1)^{\frac{75}4+15\alpha_0}}$ with an absolute constant $C$ such that $z_0=(x_0,0)$ is a regular point of  $(n,c,u)$, if
one of the following conditions holds
\ben\label{eq:condition-32}
(i)\left( \int_{Q_{1}(z_0)} n^{\frac32}(|\ln n|+1)^\frac32 + |\nabla \sqrt{c+1}|^3 + |u|^3+|p|^\frac32\right) \leq \varepsilon_2,
\een
%
%
\ben\label{eq:condition-103}
(ii)\left( \int_{Q_{1}(z_0)} n^{\frac53} + |\nabla \sqrt{c+1}|^{\frac{10}{3}} + |u|^{\frac{10}{3}}+|p|^\frac53\right) \leq \varepsilon_2'.
\een
\end{theorem}

Recall the definition of Hausdorff measure and the parabolic version:
\begin{definition} [see \cite{Tsai-2018} Chapter 6]
For a set $E\subset \mathbb{R}^{n+1}$ and $\alpha\geq 0$, $Q_r(z_0)=B_r(x_0)\times(t_0-r^2,t_0)$ for $z_0=(x_0,t_0)$. Denote by $\mathcal{P}^\alpha(E)$ its $\alpha-$dimensional parabolic Hausdorff measure, namely,
\beno
\mathcal{P}^\alpha(E)=\liminf_{\delta\rightarrow0^{+}} \left\{\sum_{j=1}^{\infty}r_j^{\alpha}: E\subset\bigcup_jQ(z_j,r_j),r_j\leq \delta\right\}.
\eeno
\end{definition}
Immediately, it follows from the above theorem that
\begin{corollary}\label{thm:3}
Taking the same assumptions as Theorem \ref{thm:regularity-fixed}, there holds $\mathcal{P}^{\frac53}(\mathcal{S}) = 0$, where $\mathcal{S}$ is the singular set at time $0$.
\end{corollary}

\begin{remark}\label{rem:3} At a fixed time, the possible singular set is described via the parabolic Hausdorff measure,
since the known a priori estimates for weak solutions are parabolic norms. The above conclusion  indicates that the concentration of cells may appear in a linear form, but not in a two-dimensional region, which is closely related to the interior singular vortices as described in \cite{DCCGK2004}.
The estimate on the Hausdorff dimension of the singular set here is weaker than the Navier-Stokes case (see \cite{CKN}). In fact, it is still unknown that whether the condition of
\ben\label{condition}
\limsup_{r \rightarrow 0} r^{-1} \int_{Q_{r}} |\nabla \sqrt{n}|^2 + |\nabla u|^2 + |\nabla^2 \sqrt{\tilde{c}}|^2 \leq \varepsilon
\een
implies the regularity, which is similar as \cite{CKN}. The main obstacle comes from the terms including $n\ln n$ on the right hand side of the local energy inequality (\ref{eq:local energy inequality}), since ``$n\ln n$" seems to be not cancelled by the term of ``$|\nabla \sqrt{n}|$" with the help of  the embedding inequality in the scaling sense.
\end{remark}

In fact, the condition of the pressure can be removed, which is stated as follows.
\begin{theorem}\label{thm:2}
Taking the same assumptions as Theorem \ref{thm:regularity-fixed}, there exists an absolute constant $\varepsilon_3$ such that$z_0=(x_0,0)$ is a regular point if,
\ben\label{eq:condition-all}
&&\limsup_{r \rightarrow 0} r^{-1}\left(\sup_{-r^2<t<0}\int_{B_{r(x_0)}} n  + |n \ln n| + |\nabla \sqrt{c+1}|^2 + |u|^2 \right)\nonumber\\
&&+\limsup_{r \rightarrow 0} r^{-1}\left( \int_{Q_{r}(z_0)} |\nabla \sqrt{n}|^2 + |\nabla u|^2 + |\nabla^2 \sqrt{c+1}|^2\right) \leq \frac{\varepsilon_3}{(\Lambda_0\Lambda_1)^{(4+\alpha_0)}}
\een
\end{theorem}


The paper is organized as follows, in Section 2, we introduce some definitions and technical lemmas, especially including the new local energy inequality. In Section 3, we  proof Theorem \ref{thm:regularity-fixed}, which is divided into four steps. Theorem {\ref{thm:lin}}, Corollary \ref{thm:3} and Theorem \ref{thm:2} are  proved in Section 4 and  Section 5, respectively.

Throughout this article, $C$ denotes an absolute constant independent of $(n,c,u)$ and may be different from line to line.

\setcounter{equation}{0}
\section{Preliminaries and some technical lemmas}

Let $(n,c,u,p)$   be a solution to the chemotaxis-Navier-Stokes equations (\ref{eq:KS}). Without loss of generality, let $z_0=(0,0)$.  Set the following scaling:
\ben\label{eq:scaling}
&&n_\lambda(x,t)=\lambda^2 n(\lambda x, \lambda^2t);~~ c_\lambda(x,t)=c(\lambda x, \lambda^2t);\nonumber\\
&&~~ u_\lambda(x,t)=\lambda u(\lambda x, \lambda^2t);~~p_\lambda(x,t)=\lambda^2 p(\lambda x, \lambda^2t)
\een
then $(n_\lambda, c_\lambda, u_\lambda,p_\lambda)$ is also a solution of (\ref{eq:KS}).

Now define some quantities which are invariant under the scaling (\ref{eq:scaling}):
\beno
&&A_u(r)=r^{-1}\|u\|^2_{L^\infty_tL^2_x(Q_r)};~~
E_u(r)=r^{-1}\|\nabla u\|^2_{L^2_tL^2_x(Q_r)};\\
&&A_c(r)=r^{-1}\|\nabla c\|^2_{L^\infty_tL^2_x(Q_r)};~~
E_c(r)=r^{-1}\|\nabla^2 c\|^2_{L^2_tL^2_x(Q_r)};\\
&&A_n(r)=r^{-1}\|\sqrt{n}\|^2_{L^\infty_tL^2_x(Q_r)};~~
E_n(r)=r^{-1}\|\nabla (\sqrt{n})\|^2_{L^2_tL^2_x(Q_r)};\\
&&C_u(r)=r^{-2}\|u\|^3_{L^3_tL^3_x(Q_r)};~~
\tilde{C}_u(r)=r^{-2}\|u-(u)_r\|^3_{L^3_tL^3_x(Q_r)};\\
&&D(r)=r^{-2}\|p\|^{\frac32}_{L^{\frac32}_tL^{\frac32}_x(Q_r)}.
\eeno

Recall a property of harmonic function.
\begin{lemma}[See \cite{Lin}]\label{mean value property} Let $f$ be a harmonic function in $B_1 \subset \mathbb{R}^3$,
for $1\leq p,q \leq \infty$, $0<r<\rho<1$ and $k\geq1$, there holds
\beno
||\nabla^k f||_{L^q(B_r)}\leq C\frac{r^\frac{3}{q}}{(\rho-r)^{\frac3p+k}}||f||_{L^p(B_\rho)}.
\eeno
\end{lemma}

\begin{lemma}[A priori estimates]\label{mean}
Under the assumptions of Theorem \ref{thm:regularity-fixed}, there holds
\beno
\nabla \sqrt{c+1} \in L^\infty L^2 \cap L^2 \dot{H}^1,
\eeno
and
\beno
\int_{\mathbb{R}^3\times (-1,0)} (\sqrt{c+1})^{-2} |\nabla \sqrt{c+1}|^4 < \infty.
\eeno
\end{lemma}
{\bf Proof.} 
 Direct calculations imply
\beno
|\nabla \sqrt{c+1}| = |\frac12(c+1)^{-\frac12}\nabla c|=\left|\frac{\sqrt{c}}{\sqrt{c+1}} \nabla \sqrt{c}\right| \leq|\nabla \sqrt{c}|,
\eeno
and
\beno
|\nabla^2 \sqrt{c+1}|\leq \frac{|\nabla \sqrt{c}|^2}{\sqrt{c+1} }+|\nabla^2 \sqrt{c}|+\frac{1}{\sqrt{c+1}}|\nabla \sqrt{c+1}||\nabla \sqrt{c}|.
\eeno
Hence we arrive at
\beno
\int_{\mathbb{R}^3 \times (-1,0)} |\nabla^2 \sqrt{c+1}|^2 \leq \int_{\mathbb{R}^3 \times (-1,0)} |\nabla^2 \sqrt{c}|^2 + \int_{\mathbb{R}^3 \times (-1,0)} \left|\frac{|\nabla \sqrt{c}|^2}{\sqrt{c+1}}\right|^2.
\eeno
Under the assumptions of Theorem \ref{thm:regularity-fixed}, there hold
\beno
\nabla \sqrt{c} \in L^\infty L^2 \cap L^2 \dot{H}^1,
\eeno
and
\beno
\int_{\mathbb{R}^3 \times (-1,0)} (\sqrt{c})^{-2} |\nabla \sqrt{c}|^4 < + \infty.
\eeno
Noting that $(\sqrt{c+1})^{-2} \leq (\sqrt{c})^{-2}$, we have
\beno
\nabla \sqrt{c+1} \in L^\infty L^2 \cap L^2 \dot{H}^1.
\eeno
The second inequality is obviously due to the relation $|\nabla \sqrt{c+1}| \leq |\nabla \sqrt{c}|$.
%

Next we establish a new local energy inequality including $\ln n$. Moreover, consider the equation of  $c+1$ instead of $c$, and we obtain some new estimates of $c+1$, which is slightly different with those in \cite{ZZ2014} and \cite{HZ2017}.
\begin{lemma}[Local energy inequality]\label{Local energy inequality}Let $\psi$ be a cut-off function, which vanishes on the parabolic boundary of ${Q_1}^t$.  Then for any $t\in (-1,0)$, the following inequality holds under the assumptions of Theorem \ref{thm:regularity-fixed}:
\ben\label{eq:local energy inequality}
&&\int_{B_1} (n \ln n \psi)(\cdot,t)dx + 2 \int_{{Q_1}^t} |\nabla \sqrt{n}|^2 \psi dxdt\nonumber\\
&&+2  \int_{B_1} (|\nabla \sqrt{\tilde{c}}|^2 \psi)(\cdot,t)dx +\frac47 \int_{{Q_1}^t} |\nabla^2 \sqrt{\tilde{c}}|^2 \psi dxdt\nonumber\\
&&+ 2 \int_{{Q_1}^t} |\nabla \sqrt{\tilde{c}}|^2 n \psi dxdt+\frac14\sum_{i,j} \int_{{Q_1}^t} (\sqrt{\tilde{c}})^{-2} (\partial_j \sqrt{\tilde{c}})^2 (\partial_i \sqrt{\tilde{c}})^2 \psi dxdt\nonumber\\
&&+112\|\tilde{c}(\cdot,-1)\|_{\infty}\int_{B_1} (|u|^2)(\cdot,t) \psi dxdt+112\|\tilde{c}(\cdot,-1)\|_{\infty}\int_{{Q_1}^t} |\nabla u|^2 \psi dxdt\\\nonumber
&&\leq \int_{{Q_1}^t} n \ln n (\partial_t \psi+ \Delta \psi)dxdt + \int_{{Q_1}^t} n \ln n u \cdot \nabla \psi dxdt\\\nonumber
&&+ \int_{{Q_1}^t} n \ln n \nabla c \cdot \nabla \psi dxdt+ \int_{{Q_1}^t} n \nabla c \cdot \nabla \psi dxdt\\\nonumber
&&+2 \int_{{Q_1}^t} |\nabla \sqrt{\tilde{c}}|^2 (\partial_t \psi + \Delta \psi)dxdt + 2 \int_{{Q_1}^t} |\nabla \sqrt{\tilde{c}}|^2 u \cdot \nabla \psi dxdt \\\nonumber
&&- \frac47\int_{{Q_1}^t} (\sqrt{\tilde{c}})^{-1} |\nabla \sqrt{\tilde{c}}|^2 \nabla \sqrt{\tilde{c}} \cdot \nabla \psi dxdt
+112\|\tilde{c}(\cdot,-1)\|_{\infty}\int_{{Q_1}^t} |u|^2 \left(\partial_t \psi + \Delta \psi\right) dxdt \\\nonumber
&&+ 112\|\tilde{c}(\cdot,-1)\|_{\infty}\int_{{Q_1}^t} |u|^2  u \cdot \nabla \psi dxdt+ 112\|\tilde{c}(\cdot,-1)\|_{\infty}\int_{{Q_1}^t} (p - \bar{p}) u \cdot \nabla \psi dxdt \\\nonumber
&&- 224\|\tilde{c}(\cdot,-1)\|_{\infty}\int_{{Q_1}^t} n \nabla \phi \cdot u \psi dxdt,
\een
where ${Q_1}^t=(-1,t)\times B_1$ and $\tilde{c} = c + 1$. 
\end{lemma}

{\bf Proof. }
Multiplying $(1+\ln n)\psi$ in $(\ref{eq:KS})_1$, integration by parts yields that
\beno
\int_{{Q_1}^t}\partial_tn (1+\ln n)\psi dxdt+\int_{{Q_1}^t}u\cdot \nabla n (1+\ln n)\psi dxdt-\int_{{Q_1}^t}\Delta n(1+\ln n)\psi dxdt\\
+\int_{{Q_1}^t}\nabla \cdot (n\nabla c) (1+\ln n)\psi dxdt\doteq T_1+\cdots+T_4=0,
\eeno
where
\beno
T_1&=&\int_{{Q_1}^t} \partial_tn{(1+\ln n)}\psi dxdt = \int_{B_1}  (n{\ln n}\psi)(\cdot,t) dx-\int_{{Q_1}^t}n\ln n\partial_t \psi dxdt;\\
T_2&=&-\int_{{Q_1}^t} u\cdot \nabla n \psi dxdt-\int_{{Q_1}^t}n(1+\ln n)u\cdot\nabla\psi dxdt=-\int_{{Q_1}^t}n\ln n u\cdot\nabla\psi dxdt;\\
T_3&=&\int_{{Q_1}^t}\frac1n \nabla n \cdot \nabla n \psi dxdt+\int_{{Q_1}^t}\nabla n \cdot \nabla \psi dxdt+\int_{{Q_1}^t} \ln n \nabla n \cdot \nabla \psi dxdt\\
&=&\int_{{Q_1}^t}\frac1n \nabla n \cdot \nabla n \psi dxdt+\int_{{Q_1}^t}\nabla n \cdot \nabla \psi dxdt-\int_{{Q_1}^t} n \frac1n \nabla n \cdot \nabla \psi dxdt-\int_{{Q_1}^t} n \ln n \Delta\psi dxdt\\
&=&4\int_{{Q_1}^t} |\nabla\sqrt{n}|^2\psi dxdt-\int_{{Q_1}^t}n\ln n\Delta\psi dxdt;\\
and\\
T_4&=&-\int_{{Q_1}^t}\nabla c\cdot \nabla n\psi dxdt-\int_{{Q_1}^t}n\nabla c\cdot\nabla\psi dxdt-\int_{{Q_1}^t}n\ln n\nabla c\cdot \nabla\psi dxdt.
\eeno
Then we have
\ben\label{ine:n 1}
&&  \int_{B_1} (n \ln n \psi)(\cdot,t)dx + 4 \int_{{Q_1}^t} |\nabla \sqrt{n}|^2 \psi dxdt\nonumber\\
&=& \int_{{Q_1}^t} n \ln n (\partial_t \psi+ \Delta \psi)dxdt + \int_{{Q_1}^t} n \ln n u \cdot \nabla \psi dxdt\nonumber\\
&&+ \int_{{Q_1}^t} n \ln n \nabla c \cdot \nabla \psi dxdt+ \int_{{Q_1}^t} \nabla n \cdot \nabla c \psi dxdt+ \int_{{Q_1}^t} n \nabla c \cdot \nabla \psi dxdt
\een

Due to $\tilde{c} = c + 1$, it follows from the equation $(\ref{eq:KS})_2$ that
\beno
\partial_t\tilde{c}+u\cdot \nabla\tilde{c}-\Delta\tilde{c}=-\tilde{c}n+n.
\eeno
Note that
\beno
\Delta \tilde{c} = 2 |\nabla \sqrt{\tilde{c}}|^2 + 2 \sqrt{\tilde{c}} \Delta \sqrt{\tilde{c}},
\eeno
and dividing $2\sqrt{\tilde{c}}$ on both sides, we get
\ben\label{eq:c-tilde}
\partial_t\sqrt{\tilde{c}}+u\cdot \nabla\sqrt{\tilde{c}}-\frac{|\nabla \sqrt{\tilde{c}}|^2}{\sqrt{\tilde{c}}}-\Delta\sqrt{\tilde{c}}= - \frac 12 \sqrt{\tilde{c}}~~n + \frac12 \frac{n}{\sqrt{\tilde{c}}}.
\een
Multiplying the above equation (\ref{eq:c-tilde}) by $-\partial_i (\partial_i \sqrt{\tilde{c}} \psi)$ and integration by parts, there holds
\beno
&&-\int_{{Q_1}^t}\partial_t\sqrt{\tilde{c}}\partial_i (\partial_i \sqrt{\tilde{c}} \psi)dxdt-\int_{{Q_1}^t}u\cdot \nabla\sqrt{\tilde{c}}\partial_i (\partial_i \sqrt{\tilde{c}} \psi)dxdt +\int_{{Q_1}^t}\frac{|\nabla \sqrt{\tilde{c}}|^2}{\sqrt{\tilde{c}}}\partial_i (\partial_i \sqrt{\tilde{c}} \psi)dxdt \\
&&+\int_{{Q_1}^t}\Delta\sqrt{\tilde{c}}\partial_i (\partial_i \sqrt{\tilde{c}} \psi)dxdt -\int_{{Q_1}^t}\frac 12 \sqrt{\tilde{c}}~~n\partial_i (\partial_i \sqrt{\tilde{c}} \psi)dxdt +  \int_{{Q_1}^t} \frac12 \frac{n}{\sqrt{\tilde{c}}}\partial_i (\partial_i \sqrt{\tilde{c}} \psi)dxdt\\
&& \doteq J_1+J_2+\cdots+J_6=0,
\eeno
where
\beno
J_1&=&\frac12\int_{{Q_1}^t}\partial_t(\partial_i\sqrt{\tilde{c}})^2\psi dxdt=\frac12\int_{B_1}(|\nabla\sqrt{\tilde{c}}|^2\psi)(\cdot,t) dx-\frac12\int_{{Q_1}^t}|\nabla\sqrt{\tilde{c}}|^2\partial_t \psi dxdt;\\
J_2&=&\int_{{Q_1}^t}\partial_iu_j\partial_j\sqrt{\tilde{c}}(\partial_i\sqrt{\tilde{c}}\psi)dxdt+\int_{{Q_1}^t}u_j\partial_{ij}\sqrt{\tilde{c}}(\partial_i\sqrt{\tilde{c}}\psi)dxdt\\
&=&\int_{{Q_1}^t} \nabla u : (\nabla\sqrt{\tilde{c}} \otimes \nabla\sqrt{\tilde{c}})\psi dxdt-\frac12 \int_{{Q_1}^t} u\cdot\nabla\psi \nabla \sqrt{\tilde{c}} \cdot \nabla \sqrt{\tilde{c}}dxdt;\\
J_3&=&-\int_{{Q_1}^t}\partial_i\big((\sqrt{\tilde{c}})^{-1}|\nabla\sqrt{\tilde{c}}|^2\big)(\partial_i\sqrt{\tilde{c}}\psi)dxdt\\
&=&\int_{{Q_1}^t}\big((\sqrt{\tilde{c}})^{-1}|\nabla\sqrt{\tilde{c}}|^2\big)\Delta\sqrt{\tilde{c}}\psi dxdt+\int_{{Q_1}^t}\big((\sqrt{\tilde{c}})^{-1}|\nabla\sqrt{\tilde{c}}|^2\big)\nabla\sqrt{\tilde{c}}\cdot\nabla\psi dxdt;\\
J_4&=&-\int_{{Q_1}^t} \partial_{ijj} \sqrt{\tilde{c}} \partial_i \sqrt{\tilde{c}} \psi dxdt=\int_{{Q_1}^t} |\nabla^2\sqrt{\tilde{c}}|^2 \psi dxdt + \int_{{Q_1}^t} \partial_{ij} \sqrt{\tilde{c}} \partial_i \sqrt{\tilde{c}} \partial_j \psi dxdt\\
&=& \int_{{Q_1}^t} |\nabla^2\sqrt{\tilde{c}}|^2 \psi dxdt - \frac12\int_{{Q_1}^t}|\nabla\sqrt{\tilde{c}}|^2\Delta \psi dxdt;\\
J_5&=&\frac12\int_{{Q_1}^t}\sqrt{\tilde{c}}n(-\partial_i(\partial_i\sqrt{\tilde{c}}\psi))dxdt\\
&=&\frac12\int_{{Q_1}^t}|\nabla\sqrt{\tilde{c}}|^2n\psi dxdt +\frac12\int_{{Q_1}^t}\sqrt{\tilde{c}}\nabla n\cdot\nabla\sqrt{\tilde{c}} \psi dxdt\\
&=&\frac12\int_{{Q_1}^t}|\nabla\sqrt{\tilde{c}}|^2n\psi dxdt +\frac14\int_{{Q_1}^t}\nabla c \cdot \nabla n\psi dxdt;\\
\eeno
and
\beno
J_6&=&-\frac12 \int_{{Q_1}^t}\nabla\left( \frac{n}{\sqrt{\tilde{c}}}\right) \cdot \nabla \sqrt{\tilde{c}} \psi dxdt
\eeno
Combining all of them, we have
\ben\label{ine:c 1}  \nonumber
&&\frac12   \int_{B_1} (|\nabla \sqrt{\tilde{c}}|^2 \psi)(\cdot,t)dx + \int_{{Q_1}^t} |\nabla^2 \sqrt{\tilde{c}}|^2 \psi dxdt\nonumber\\
&=& \frac12 \int_{{Q_1}^t} |\nabla \sqrt{\tilde{c}}|^2 (\partial_t \psi + \Delta \psi)dxdt + \frac 12 \int_{{Q_1}^t} |\nabla \sqrt{\tilde{c}}|^2 u \cdot \nabla \psi dxdt\nonumber \\
  \nonumber
&&- \int_{{Q_1}^t} \nabla \sqrt{\tilde{c}} \cdot \nabla u \cdot \nabla \sqrt{\tilde{c}} \psi dxdt - \int_{{Q_1}^t} (\sqrt{\tilde{c}})^{-1} |\nabla \sqrt{\tilde{c}}|^2 \Delta \sqrt{\tilde{c}} \psi dxdt\\  \nonumber
&&- \int_{{Q_1}^t} (\sqrt{\tilde{c}})^{-1} |\nabla \sqrt{\tilde{c}}|^2 \nabla \sqrt{\tilde{c}} \cdot \nabla \psi dxdt - \frac12 \int_{{Q_1}^t} |\nabla \sqrt{\tilde{c}}|^2 n \psi dxdt\\
&&- \frac14 \int_{{Q_1}^t} \nabla n \cdot \nabla c \psi dxdt+ \frac12 \int_{{Q_1}^t} \nabla \left(\frac{n}{\sqrt{\tilde{c}}}\right) \cdot \nabla \sqrt{\tilde{c}} \psi dxdt.
\een
We remark here the bad terms are those integrals without $\nabla\psi$.
One bad term of all the above terms is $I\doteq-\int_{{Q_1}^t} (\sqrt{\tilde{c}})^{-1} |\nabla \sqrt{\tilde{c}}|^2 \Delta \sqrt{\tilde{c}} \psi$, and next we estimate it in details. Firstly, integration by parts yields that
\beno
I &=& - \int_{{Q_1}^t} (\sqrt{\tilde{c}})^{-1} (\partial_j \sqrt{\tilde{c}})^2 \partial_{ii} \sqrt{\tilde{c}} \psi dxdt \\
&=& - \int_{{Q_1}^t} (\sqrt{\tilde{c}})^{-2} (\partial_j \sqrt{\tilde{c}})^2 (\partial_i \sqrt{\tilde{c}})^2 \psi dxdt + 2 \int_{{Q_1}^t} (\sqrt{\tilde{c}})^{-1} \partial_{ij} \sqrt{\tilde{c}} \partial_i \sqrt{\tilde{c}} \partial_j \sqrt{\tilde{c}} \psi dxdt \\
&&+ \int_{{Q_1}^t} (\sqrt{\tilde{c}})^{-1} |\nabla \sqrt{\tilde{c}}|^2 \nabla \sqrt{\tilde{c}} \cdot \nabla \psi dxdt \\
&=& - \sum_{i,j} \int_{{Q_1}^t} (\sqrt{\tilde{c}})^{-2} (\partial_j \sqrt{\tilde{c}})^2 (\partial_i \sqrt{\tilde{c}})^2 \psi dxdt + 2 \sum_{i = j} \int_{{Q_1}^t} (\sqrt{\tilde{c}})^{-1} \partial_{ij} \sqrt{\tilde{c}} \partial_i \sqrt{\tilde{c}} \partial_j \sqrt{\tilde{c}} \psi dxdt \\
&&+ 2 \sum_{i \neq j} \int_{{Q_1}^t} (\sqrt{\tilde{c}})^{-1} \partial_{ij} \sqrt{\tilde{c}} \partial_i \sqrt{\tilde{c}} \partial_j \sqrt{\tilde{c}} \psi dxdt  + \int_{{Q_1}^t} (\sqrt{\tilde{c}})^{-1} |\nabla \sqrt{\tilde{c}}|^2 \nabla \sqrt{\tilde{c}} \cdot \nabla \psi dxdt.
\eeno
Noting that
\beno
\sum_{i = j} \int_{{Q_1}^t} (\sqrt{\tilde{c}})^{-1} \partial_{ij} \sqrt{\tilde{c}} \partial_i \sqrt{\tilde{c}} \partial_j \sqrt{\tilde{c}} \psi dxdt = - I - \sum_{i \neq j} \int_{{Q_1}^t} (\sqrt{\tilde{c}})^{-1} \partial_{ii} \sqrt{\tilde{c}} \partial_j \sqrt{\tilde{c}} \partial_j \sqrt{\tilde{c}} \psi dxdt,
\eeno
we have
\beno
I &=& - \frac 13 \sum_{i,j} \int_{{Q_1}^t} (\sqrt{\tilde{c}})^{-2} (\partial_j \sqrt{\tilde{c}})^2 (\partial_i \sqrt{\tilde{c}})^2 \psi dxdt - \frac 23 \sum_{i \neq j} \int_{{Q_1}^t} (\sqrt{\tilde{c}})^{-1} \partial_{ii} \sqrt{\tilde{c}} \partial_j \sqrt{\tilde{c}} \partial_j \sqrt{\tilde{c}} \psi dxdt \\
&&+ \frac 23 \sum_{i \neq j} \int_{{Q_1}^t} (\sqrt{\tilde{c}})^{-1} \partial_{ij} \sqrt{\tilde{c}} \partial_i \sqrt{\tilde{c}} \partial_j \sqrt{\tilde{c}} \psi dxdt+ \frac13 \int_{{Q_1}^t} (\sqrt{\tilde{c}})^{-1} |\nabla \sqrt{\tilde{c}}|^2 \nabla \sqrt{\tilde{c}} \cdot \nabla \psi dxdt.
\eeno
Secondly, using Young inequality of $ab\leq \epsilon{a^2}+\frac{b^2}{4\epsilon}$, it follows that
\beno
\frac 23 \sum_{i \neq j} \int_{{Q_1}^t} (\sqrt{\tilde{c}})^{-1} \partial_{ij} \sqrt{\tilde{c}} \partial_i \sqrt{\tilde{c}} \partial_j \sqrt{\tilde{c}} \psi dxdt &\leq& \frac 13 \sum_{i \neq j} \int_{{Q_1}^t} (\sqrt{\tilde{c}})^{-2} |\partial_i \sqrt{\tilde{c}}|^2 |\partial_j \sqrt{\tilde{c}}|^2 \psi dxdt \\
&&+ \frac 13 \sum_{i \neq j} \int_{{Q_1}^t} |\partial_{ij} \sqrt{\tilde{c}}|^2 \psi dxdt,
\eeno
and
\beno
&&- \frac 23 \sum_{i \neq j} \int_{{Q_1}^t} (\sqrt{\tilde{c}})^{-1} \partial_{ii} \sqrt{\tilde{c}} \partial_j \sqrt{\tilde{c}} \partial_j \sqrt{\tilde{c}} \psi\\
&=&- \frac 23 \sum_{i\neq j} \int_{{Q_1}^t} (\sqrt{\tilde{c}})^{-2} (\partial_j \sqrt{\tilde{c}})^2 (\partial_i \sqrt{\tilde{c}})^2 \psi +\frac 43 \sum_{i \neq j} \int_{{Q_1}^t} (\sqrt{\tilde{c}})^{-1} \partial_{ij} \sqrt{\tilde{c}} \partial_i \sqrt{\tilde{c}} \partial_j \sqrt{\tilde{c}} \psi\\
&&+ \frac23 \sum_{i\neq j}  \int_{{Q_1}^t} (\sqrt{\tilde{c}})^{-1} |\nabla_j \sqrt{\tilde{c}}|^2 \nabla_i \sqrt{\tilde{c}}  \nabla_i \psi\\
\\ &\leq& \frac 23  \sum_{i\neq j} \int_{{Q_1}^t} |\partial_{ij} \sqrt{\tilde{c}}|^2 \psi+\frac23 \int_{{Q_1}^t} (\sqrt{\tilde{c}})^{-1} |\nabla \sqrt{\tilde{c}}|^2 \nabla \sqrt{\tilde{c}} \cdot \nabla \psi
\eeno
Then we have
\ben\label{eq:I}
I &\leq& - \frac 1{3} \sum_{i=j} \int_{{Q_1}^t} (\sqrt{\tilde{c}})^{-2} (\partial_j \sqrt{\tilde{c}})^2 (\partial_i \sqrt{\tilde{c}})^2 \psi +  \sum_{i\neq j} \int_{{Q_1}^t} |\partial_{ij} \sqrt{\tilde{c}}|^2 \psi \nonumber\\
&&+ \int_{{Q_1}^t} (\sqrt{\tilde{c}})^{-1} |\nabla \sqrt{\tilde{c}}|^2 \nabla \sqrt{\tilde{c}} \cdot \nabla \psi
\een
Submitting it to (\ref{ine:c 1}), we get
\ben\label{ine:c 1'}  \nonumber
&&\frac12   \int_{B_1} (|\nabla \sqrt{\tilde{c}}|^2 \psi)(\cdot,t) + \frac13\int_{{Q_1}^t} |\triangle \sqrt{\tilde{c}}|^2 \psi \nonumber\\
&&+ \frac12 \int_{{Q_1}^t} |\nabla \sqrt{\tilde{c}}|^2 n \psi+\frac 1{3} \sum_{i=j} \int_{{Q_1}^t} (\sqrt{\tilde{c}})^{-2} (\partial_j \sqrt{\tilde{c}})^2 (\partial_i \sqrt{\tilde{c}})^2 \psi\nonumber\\
&&\leq \frac12 \int_{{Q_1}^t} |\nabla \sqrt{\tilde{c}}|^2 (\partial_t \psi + \Delta \psi) + \frac 12 \int_{{Q_1}^t} |\nabla \sqrt{\tilde{c}}|^2 u \cdot \nabla \psi\nonumber \\
  \nonumber
&&- \int_{{Q_1}^t} \nabla \sqrt{\tilde{c}} \cdot \nabla u \cdot \nabla \sqrt{\tilde{c}} \psi  \\
&&- \frac14 \int_{{Q_1}^t} \nabla n \cdot \nabla c \psi + \frac12 \int_{{Q_1}^t} \nabla \left(\frac{n}{\sqrt{\tilde{c}}}\right) \cdot \nabla \sqrt{\tilde{c}} \psi.
\een
Note that
\beno
\int_{{Q_1}^t} (\sqrt{\tilde{c}})^{-2} |\nabla \sqrt{\tilde{c}}|^4 \psi&=& \sum_{i,j} \int_{{Q_1}^t} (\sqrt{\tilde{c}})^{-2} (\partial_j \sqrt{\tilde{c}})^2 (\partial_i \sqrt{\tilde{c}})^2 \psi\\
&\leq&  3\sum_{i=j} \int_{{Q_1}^t} (\sqrt{\tilde{c}})^{-2} (\partial_j \sqrt{\tilde{c}})^2 (\partial_i \sqrt{\tilde{c}})^2 \psi
\eeno
and we have
\ben\label{eq: estimate of I}
I&=&- \sum_{i,j}\int_{{Q_1}^t} (\sqrt{\tilde{c}})^{-1} (\partial_j \sqrt{\tilde{c}})^2 \partial_{ii} \sqrt{\tilde{c}} \psi\nonumber \\
&\leq &\frac32\int_{{Q_1}^t} |\triangle \sqrt{\tilde{c}}|^2 \psi+\frac 1{2} \sum_{i=j} \int_{{Q_1}^t} (\sqrt{\tilde{c}})^{-2} (\partial_j \sqrt{\tilde{c}})^2 (\partial_i \sqrt{\tilde{c}})^2 \psi.
\een
Then with the help of (\ref{eq: estimate of I}), by  $(\ref{ine:c 1})+6\times (\ref{ine:c 1'})$ we get
\ben\label{ine:c 1''}
&&\frac72   \int_{B_1} (|\nabla \sqrt{\tilde{c}}|^2 \psi)(\cdot,t)dx + \int_{{Q_1}^t} |\nabla^2 \sqrt{\tilde{c}}|^2 \psi dxdt\nonumber\\
&&+ \frac72  \int_{{Q_1}^t} |\nabla \sqrt{\tilde{c}}|^2 n \psi dxdt+\frac12\int_{{Q_1}^t} (\sqrt{\tilde{c}})^{-2} |\nabla \sqrt{\tilde{c}}|^4 \psi dxdt
\nonumber\\
&&\leq \frac72 \int_{{Q_1}^t} |\nabla \sqrt{\tilde{c}}|^2 (\partial_t \psi + \Delta \psi)dxdt + \frac 72 \int_{{Q_1}^t} |\nabla \sqrt{\tilde{c}}|^2 u \cdot \nabla \psi dxdt\nonumber \\
  \nonumber
&&- 7\int_{{Q_1}^t} \nabla \sqrt{\tilde{c}} \cdot \nabla u \cdot \nabla \sqrt{\tilde{c}} \psi dxdt
- \int_{{Q_1}^t} (\sqrt{\tilde{c}})^{-1} |\nabla \sqrt{\tilde{c}}|^2 \nabla \sqrt{\tilde{c}} \cdot \nabla \psi dxdt \\
&&- \frac74 \int_{{Q_1}^t} \nabla n \cdot \nabla c \psi dxdt+ \frac72 \int_{{Q_1}^t} \nabla \left(\frac{n}{\sqrt{\tilde{c}}}\right) \cdot \nabla \sqrt{\tilde{c}} \psi dxdt,
\een
where we neglected the term of $\frac12\int_{{Q_1}^t} |\triangle \sqrt{\tilde{c}}|^2 \psi dxdt.$
%

Moreover,
using Young's inequality
\ben\label{ine:g}
\left|\int_{{Q_1}^t} \nabla \sqrt{\tilde{c}} \cdot \nabla u \cdot \nabla \sqrt{\tilde{c}} \psi\right| &\leq& \frac1{112} \int_{{Q_1}^t} (\sqrt{\tilde{c}})^{-2} |\nabla \sqrt{\tilde{c}}|^4\psi + 28 \int_{{Q_1}^t} (\sqrt{\tilde{c}})^2|\nabla u|^2 \psi\nonumber\\
&\leq & \frac1{112} \int_{{Q_1}^t} (\sqrt{\tilde{c}})^{-2} |\nabla \sqrt{\tilde{c}}|^4\psi + 28\|\tilde{c}\|_{\infty} \int_{{Q_1}^t} |\nabla u|^2 \psi\nonumber\\
\een
and
\ben\label{eq:III}
&& \frac72\int_{{Q_1}^t} \nabla \left(\frac{n}{\sqrt{\tilde{c}}}\right) \cdot \nabla \sqrt{\tilde{c}} \psi\nonumber\\
  &=& -\frac72\int_{{Q_1}^t} (\sqrt{\tilde{c}})^{-2} n |\nabla \sqrt{\tilde{c}}|^2 \psi + 7\int_{{Q_1}^t} (\sqrt{\tilde{c}})^{-1} \sqrt{n} \nabla \sqrt{n} \cdot \nabla \sqrt{\tilde{c}} \psi\nonumber \\
&\leq& -\frac72\int_{{Q_1}^t} (\sqrt{\tilde{c}})^{-2} n |\nabla \sqrt{\tilde{c}}|^2 \psi + \frac72\int_{{Q_1}^t} (\sqrt{\tilde{c}})^{-2} n |\nabla \sqrt{\tilde{c}}|^2 \psi + \frac72\int_{{Q_1}^t} |\nabla \sqrt{n}|^2 \psi\nonumber\\
&\leq& \frac72\int_{{Q_1}^t} |\nabla \sqrt{n}|^2 \psi,
\een

Substitute these estimates (\ref{ine:g})-(\ref{eq:III}) to the inequality of (\ref{ine:c 1''}), and we get

\ben\label{ine:c 2}
&&\frac72   \int_{B_1} (|\nabla \sqrt{\tilde{c}}|^2 \psi)(\cdot,t)dx + \int_{{Q_1}^t} |\nabla^2 \sqrt{\tilde{c}}|^2 \psi dxdt\nonumber\\
&&+ \frac72  \int_{{Q_1}^t} |\nabla \sqrt{\tilde{c}}|^2 n \psi+\frac7{16}\sum_{i,j} \int_{{Q_1}^t} (\sqrt{\tilde{c}})^{-2} (\partial_j \sqrt{\tilde{c}})^2 (\partial_i \sqrt{\tilde{c}})^2 \psi\nonumber\\
&&\leq \frac72 \int_{{Q_1}^t} |\nabla \sqrt{\tilde{c}}|^2 (\partial_t \psi + \Delta \psi)dxdt + \frac 72 \int_{{Q_1}^t} |\nabla \sqrt{\tilde{c}}|^2 u \cdot \nabla \psi dxdt\nonumber \\
  \nonumber
&&- \frac74 \int_{{Q_1}^t} \nabla n \cdot \nabla c \psi dxdt
- \int_{{Q_1}^t} (\sqrt{\tilde{c}})^{-1} |\nabla \sqrt{\tilde{c}}|^2 \nabla \sqrt{\tilde{c}} \cdot \nabla \psi dxdt \\
&&+196\|\tilde{c}\|_{\infty} \int_{{Q_1}^t} |\nabla u|^2 \psi+ \frac72\int_{{Q_1}^t} |\nabla \sqrt{n}|^2 \psi dxdt.
\een

Recall the local estimate of $n$ in (\ref{ine:n 1}), by taking
$(\ref{ine:n 1}) + \frac47\times(\ref{ine:c 2})$, we arrive at
\ben\label{ine:energy n c}\nonumber
&&\int_{B_1} (n \ln n \psi)(\cdot,t)dx + 2 \int_{{Q_1}^t} |\nabla \sqrt{n}|^2 \psi dxdt\nonumber\\
&&+2  \int_{B_1} (|\nabla \sqrt{\tilde{c}}|^2 \psi)(\cdot,t)dx +\frac47 \int_{{Q_1}^t} |\nabla^2 \sqrt{\tilde{c}}|^2 \psi dxdt\nonumber\\
&&+ 2 \int_{{Q_1}^t} |\nabla \sqrt{\tilde{c}}|^2 n \psi+\frac14\sum_{i,j} \int_{{Q_1}^t} (\sqrt{\tilde{c}})^{-2} (\partial_j \sqrt{\tilde{c}})^2 (\partial_i \sqrt{\tilde{c}})^2 \psi \nonumber\\
&&\leq \int_{{Q_1}^t} n \ln n (\partial_t \psi+ \Delta \psi)dxdt + \int_{{Q_1}^t} n \ln n u \cdot \nabla \psi dxdt\\\nonumber
&&+ \int_{{Q_1}^t} n \ln n \nabla c \cdot \nabla \psi dxdt+ \int_{{Q_1}^t} n \nabla c \cdot \nabla \psi dxdt\\\nonumber
&&+2 \int_{{Q_1}^t} |\nabla \sqrt{\tilde{c}}|^2 (\partial_t \psi + \Delta \psi)dxdt + 2 \int_{{Q_1}^t} |\nabla \sqrt{\tilde{c}}|^2 u \cdot \nabla \psi dxdt \\\nonumber
&&- \frac47\int_{{Q_1}^t} (\sqrt{\tilde{c}})^{-1} |\nabla \sqrt{\tilde{c}}|^2 \nabla \sqrt{\tilde{c}} \cdot \nabla \psi dxdt+112\|\tilde{c}\|_{\infty} \int_{{Q_1}^t} |\nabla u|^2 \psi.
\een

Multiplying the third equation of (\ref{eq:KS}) by $2u \psi$ and integration by parts, we have
\ben\label{ine:energy u}\nonumber
\int_{B_1} (|u|^2)(\cdot,t) \psi + 2\int_{{Q_1}^t} |\nabla u|^2 \psi &&\leq \int_{{Q_1}^t} |u|^2 \left(\partial_t \psi + \Delta \psi\right) + \int_{{Q_1}^t} |u|^2  u \cdot \nabla \psi\\
&& + \int_{{Q_1}^t} (p - \bar{p}) u \cdot \nabla \psi - 2\int_{{Q_1}^t} n \nabla \phi \cdot u \psi
\een
where $\bar{p}$ is independent of the space variable $x$.
Taking  $(\ref{ine:energy n c}) +112\|\tilde{c}(\cdot,-1)\|_{\infty}\times (\ref{ine:energy u})$, which yields (\ref{eq:local energy inequality}) via maximum principle of $c$'s equation. The proof is complete.

\section{proof of Theorem {\ref{thm:regularity-fixed}}}

Without loss of generality, let $z_0=(0,0)$.
Before we begin the proof, for $r_k=2^{-k}$ with $k\in { \mathbb{N}}$, we introduce the following backward heat kernel $\Psi_{n}$ as in \cite{CKN}:
\beno
\Psi_{n}(x,t) = \frac{1}{(r_{n}^2 - t )^\frac 32} \exp(- \frac{|x |^2}{4(r_{n}^2 - t )}),
\eeno
where $(x,t) \in \mathbb{R}^3 \times (-\infty,r_{n}^2)$. 
%
Take a suitable cut-off function $\xi(x,t)$ in $ Q_{r_3}$, which satisfies
\begin{eqnarray}\label{cut off}
 \xi(x,t)=\left\{
    \begin{array}{llll}
    \displaystyle 1,\quad {\rm in}\quad Q_{r_4} \\
    \displaystyle 0,\quad {\rm in}\quad Q_{r_3}^{c}.\\
    \end{array}
 \right.
\end{eqnarray}
It is easy to check the following properties of $\phi_n = \Psi_n \xi$.
\begin{proposition}\label{prop:1}
There exist two absolute constants  $C_1$ and $C_2$ such that\\
(i). $C_1 r_{n}^{-3} \leq \phi_n(x,t) \leq C_2 r_{n}^{-3}$ on $Q_{r_{n}}$ for $n\geq 2$;\\
(ii). $\phi_n(x,t) \leq C_2 r_k^{-3}$ for $(x,t) \in Q_{r_k} \setminus Q_{r_{k+1}}$, \quad $1<k\leq n$;\\
(iii).$|\nabla \phi_n(x,t) |\leq C_2r_n^{-4}$ in $Q_{r_n}$,\quad $n\geq 2$;\\
(iv). $|\nabla \phi_n(x,t)| \leq C_2 r_k^{-4}$ on $Q_{r_{k-1}} \setminus Q_{r_{k}}$, \quad $1<k\leq n$;\\
(v). $\left|\partial_t \phi + \Delta \phi\right| \leq C$ on $Q_{r_3}$;\\
(vi). $\partial_t \phi + \Delta \phi=0$ on $Q_{r_4}$.
\end{proposition}

{\bf Proof of Theorem {\ref{thm:regularity-fixed}}.}

In the following proof, by mathematical induction
 we are aimed to prove the following inequality
\ben\label{ine:induction}
&&r_k^{-3}\sup_{-r_k^2<t<0} \int_{B_{r_k}} n  + |n \ln n| + |\nabla \sqrt{\tilde{c}}|^2 + |u|^2 \nonumber\\
&&+ r_k^{-3} \int_{Q_{r_k}} |\nabla \sqrt{n}|^2 + |\nabla^2 \sqrt{\tilde{c}}|^2 + |\nabla u|^2 \leq C_0 \varepsilon_0^\frac12,
\een
for any $k\geq 1,$ where $C_0>1$ is an absolute constant and $ \varepsilon_0=\frac{\varepsilon_1}{(\Lambda_0\Lambda_1)^{(4+\alpha_0)}}$.
Obviously, (\ref{eq:condition-fixed}) implies that (\ref{ine:induction}) holds for $k = 1$. Assume that (\ref{ine:induction}) holds for the case of  $k=1,2,\cdots, N$. Next we prove (\ref{ine:induction}) the case of $k=N+1$.

{\bf Step 1: Estimates from the local energy inequality.}

Taking $\psi=\phi_{N+1}$ as a test function in the local energy inequality (\ref{eq:local energy inequality}), we have
\beno \nonumber
&&\int_{B_{r_{N+1}}} (n \ln n \psi)(\cdot,t)  + r_{N+1}^{-3}\int_{B_{r_{N+1}}} |\nabla \sqrt{\tilde{c}}|^2(\cdot,t) + \Lambda_0r_{N+1}^{-3}\int_{B_{r_{N+1}}} |u|^2(\cdot,t)  \\
&&+ r_{N+1}^{-3}\int_{{Q_{r_{N+1}}}^t} |\nabla \sqrt{n}|^2 + \Lambda_0 r_{N+1}^{-3}\int_{{Q_{r_{N+1}}}^t} |\nabla u|^2\nonumber\\
&&+ r_{N+1}^{-3}\int_{{Q_{r_{N+1}}}^t} |\nabla^2 \sqrt{\tilde{c}}|^2 + r_{N+1}^{-3}\int_{{Q_{r_{N+1}}}^t} (\sqrt{\tilde{c}})^{-2} |\nabla \sqrt{\tilde{c}}|^4 +  r_{N+1}^{-3}\int_{{Q_{r_{N+1}}}^t} |\nabla \sqrt{\tilde{c}}|^2 n  \nonumber\\ \nonumber
&\leq& C \int_{Q_{r_3}}\left| n \ln n (\partial_t \psi + \Delta \psi)\right| + C\int_{Q_{r_3}}\left| n \ln n u \cdot \nabla \psi \right|+ C \int_{Q_{r_3}}\left| n \ln n \nabla c \cdot \nabla \psi \right|\\ \nonumber
&&+ C \int_{Q_{r_3}}\left| n \nabla c \cdot \nabla \psi\right| + C \int_{Q_{r_3}} \left||\nabla \sqrt{\tilde{c}}|^2 (\partial_t \psi + \Delta \psi)\right| + C \int_{Q_{r_3}}\left| |\nabla \sqrt{\tilde{c}}|^2 u \cdot \nabla \psi \right|\\ \nonumber
&&+ C \int_{Q_{r_3}} (\sqrt{\tilde{c}})^{-1} |\nabla \sqrt{\tilde{c}}|^2 |\nabla \sqrt{\tilde{c}} \cdot \nabla \psi| + C \Lambda_0\int_{Q_{r_3}} |u|^2 \left|\left(\partial_t \psi + \Delta \psi\right) \right|\\
&&+ C \Lambda_0\int_{Q_{r_3}} |u|^2 \left| u \cdot \nabla \psi \right|+ C\Lambda_0 \left|\int_{-r_3^2}^t\int_{B_{r_3}} (p - \bar{p}) u \cdot \nabla \psi \right|+ C \Lambda_0\int_{Q_{r_3}}\left| n \nabla \phi \cdot u \psi\right|\nonumber \\
&&:= I_1 + I_2 + \cdots + I_{11}.
\eeno

For the term $I_1$, $I_2$ and $I_3$, we need to deal with the part $n \ln n$. For $1\leq k\leq N$, by the embedding inequality
\beno
||n^{\frac12}||_{L^\frac{10}3_{t,x}(Q_{r_{k}})} \leq C||n^{\frac12}||^{\frac25}_{L^\infty_tL^2_x(Q_{r_{k}})}||\nabla n^{\frac12}||^{\frac35}_{L^2_tL^2_x(Q_{r_{k}})}+C||n^{\frac12}||_{L^\infty_tL^2_x(Q_{r_{k}})}
\eeno
and by (\ref{ine:induction}) we have
\ben\label{ine:n 5 3}
r_{k}^{-3} \|n\|_{L^\frac53(Q_{r_{k}})} \leq C_0 C\varepsilon_0^\frac12.
\een
Decomposing $Q_{r_{k}}$ into $Q_{r_{k}} \cap \{n \leq A\}$ and $Q_{r_{k}} \cap \{n>A\}$ with a constant $A$, we arrive at
\beno
\int_{Q_{r_{k}}} |n \ln n|^{\frac53-\delta} &=& \int_{Q_{r_{k}} \cap \{n(x) \leq A\}} |n \ln n|^{\frac53-\delta} + \int_{Q_{r_{k}} \cap \{ n(x) > A\}} |n \ln n|^{\frac53-\delta}.
\eeno
Since
\ben\label{n_0}
\lim_{n \rightarrow 0} n^\frac{\delta}{\frac53-\delta} \ln n = 0,
\een
we know that  for $0<\delta<1/3$, $n^\delta |\ln n|^{\frac53-\delta} \leq C(\delta,A)$ in the domain $\{x:n(x) \leq A\}$.

On the other hand,
\ben\label{n_1}
\lim_{n \rightarrow \infty} n^{-\delta}|\ln n|^{\frac53-\delta} = 0.
\een
Thus $n^{-\delta} |\ln n|^{\frac53-\delta} \leq C(\delta,A)$ in the domain $\{x:n(x) > A\}$ for a large constant $A$.
Fixed $A$ and $\delta$, for example one can choose $A=100$ and $\delta=\frac16.$
From (\ref{n_0}) and (\ref{n_1}), we have
\beno\label{ine:estimate n ln n}
\int_{Q_{r_{k}}} |n \ln n|^{\frac32} \leq C \int_{Q_{r_{k}} \cap \{ n(x) \leq 100\}} |n|^{\frac43} + C \int_{Q_{r_{k}} \cap \{ n(x) > 100\}} |n|^\frac53,
\eeno
which is controlled by
\ben\label{ine:n ln n}
\int_{Q_{r_{k}}} |n \ln n|^{\frac32}&\leq &C \int_{Q_{r_{k}}} |n|^{\frac43} + r_{k}^5C C_0^\frac53 \varepsilon_0^\frac56 \nonumber\\
&\leq& C r_{k}^5 C_0^{\frac43} \varepsilon_0^{\frac23} + C r_{k}^5 C_0^\frac53 \varepsilon_0^\frac56\leq  C r_{k}^5 C_0^\frac53 \varepsilon_0^\frac23,
\een
due to  (\ref{ine:induction}) and H\"{o}lder inequality.

{\bf \underline{Estimate of $I_1$}. }
Noting that Proposition \ref{prop:1} and (\ref{eq:condition-fixed}), we have
\beno
I_1 \leq  C \int_{Q_{r_3}} |n \ln n| \leq C \varepsilon_0.
\eeno

{\bf \underline{Estimate of $I_2$}. } Using Proposition \ref{prop:1}, (\ref{ine:induction}) and (\ref{ine:n ln n}), we have
\beno
I_2 &\leq & C\sum_{k = 1}^N \int_{Q_{r_k} \setminus Q_{r_{k+1}}} \left|n \ln n u \cdot \nabla \psi\right| + C\int_{Q_{r_{N+1}}} \left|n \ln n u \cdot \nabla \psi\right|\\
&\leq& \sum_{k=1}^N C r_k^{-4} \|n \ln n\|_{L^{\frac32}(Q_{r_k})} \|u\|_{L^\frac{10}3(Q_{r_k})} r_k^{\frac16} + C r_{N+1}^{-4} \|n \ln n\|_{L^{\frac32}(Q_{r_N})} \|u\|_{L^\frac{10}3(Q_{r_N})} r_{N+1}^{\frac16} \\
&\leq& C C_0^\frac{29}{18} \sum_{k=1}^N r_k^{-4} r_k^{\frac16} r_k^\frac{10}{3} \varepsilon_0^\frac{4}{9} r_k^{\frac32} \varepsilon_0^\frac14 +C C_0^\frac{29}{18} r_{N+1}^{-4} r_{N+1}^{\frac16} r_N^\frac{10}{3} r_N^\frac32 \varepsilon_0^\frac94 \varepsilon_0^\frac14 \\
&\leq& C C_0^\frac{29}{18} \varepsilon_0^\frac{25}{36}.
\eeno

{\bf \underline{Estimate of $I_3$}. } Since $u$ is similar as $\nabla \sqrt{\tilde{c}}$, using Proposition \ref{prop:1}, (\ref{ine:induction}) and (\ref{ine:n ln n}) agian, we have
\beno
I_3 \leq C C_0^\frac{29}{18} \varepsilon_0^\frac{25}{36}\sqrt{\Lambda_0}.
\eeno

{\bf \underline{Estimate of $I_4$}. }
Using Proposition \ref{prop:1}, (\ref{eq:condition-fixed}) and (\ref{ine:n 5 3}), we have
\beno
I_4 &\leq & C\sum_{k=1}^N \int_{Q_{r_k} \setminus Q_{r_{k+1}}} \left|n \nabla c \cdot \nabla \psi\right| + C\int_{Q_{r_{N+1}}}\left| n \nabla c \cdot \nabla \psi \right| \\
&\leq& C\sqrt{\Lambda_0} \sum_{k=1}^N r_k^{-4} \|n\|_{L^\frac53(Q_{r_k})} \|\nabla \sqrt{\tilde{c}}\|_{L^\frac{10}3(Q_{r_k})} r_k^\frac12 + C r_{N+1}^{-4} \|n\|_{L^\frac53(Q_{r_N})} \|\nabla \sqrt{\tilde{c}}\|_{L^\frac{10}3(Q_{r_N})} r_{N+1}^\frac12 \\
&\leq& C\sqrt{\Lambda_0} C_0^\frac32 \sum_{k=1}^N r_k^{-4} r_k^3 \varepsilon_0^\frac12 r_k^\frac32 \varepsilon_0^\frac14 r_k^\frac12 + C C_0^\frac32 r_{N+1}^{-4} r_N^3 \varepsilon_0^\frac12 r_N^\frac32 \varepsilon_0^\frac14 r_k^\frac12 \\
&\leq& C \sqrt{\Lambda_0}C_0^\frac32 \varepsilon_0^\frac34.
\eeno

{\bf \underline{Estimate of $I_5$}. }
Noting that Proposition \ref{prop:1} and (\ref{eq:condition-fixed}), we have
\beno
I_5 \leq C \int_{Q_{r_3}} |\nabla \sqrt{\tilde{c}}|^2 \leq C \varepsilon_0.
\eeno

{\bf \underline{Estimate of $I_6$}. }
Similar as the estimate of $I_4$, by Proposition \ref{prop:1}, (\ref{ine:induction}) and embedding inequality, we have
\beno
I_6 &\leq & C\sum_{k=1}^N \int_{Q_{r_k} \setminus Q_{r_{k+1}}} |\nabla \sqrt{\tilde{c}}|^2 \left|u \cdot \nabla \psi\right| + C\int_{Q_{r_{N+1}}} |\nabla \sqrt{\tilde{c}}|^2 \left|u \cdot \nabla \psi \right|\\
&\leq& C \sum_{k=1}^N r_k^{-4} \|\nabla \sqrt{\tilde{c}}\|_{L^\frac{10}3(Q_{r_k})}^2 \|u\|_{L^\frac{10}3(Q_{r_k})} r_k^\frac12 + C r_{N+1}^{-4} \|\nabla \sqrt{\tilde{c}}\|_{L^\frac{10}3(Q_{r_N})}^2 \|u\|_{L^\frac{10}3(Q_{r_N})} r_{N+1}^\frac12 \\
&\leq& C C_0^\frac32 \sum_{k=1}^N r_k^{-4} r_k^3 \varepsilon_0^\frac12 r_k^\frac32 \varepsilon_0^\frac14 r_k^\frac12 + C C_0^\frac32 r_{N+1}^{-4} r_N^3 \varepsilon_0^\frac12 r_N^\frac32 \varepsilon_0^\frac14 r_{N+1}^\frac12 \\
&\leq& C C_0^\frac32 \varepsilon_0^\frac34.
\eeno

{\bf \underline{Estimate of $I_7$}. }
Noting that $\tilde c\geq 1$, we have
\beno
I_7 &\leq & C\sum_{k=1}^N \int_{Q_{r_k} \setminus Q_{r_{k+1}}} (\sqrt{\tilde{c}})^{-1} |\nabla \sqrt{\tilde{c}}|^2 |\nabla \sqrt{\tilde{c}} \cdot \nabla \psi |+ C\int_{Q_{r_{N+1}}} (\sqrt{\tilde{c}})^{-1} |\nabla \sqrt{\tilde{c}}|^2 |\nabla \sqrt{\tilde{c}} \cdot \nabla \psi |\\
&\leq& C \sum_{k=1}^N r_k^{-4} \|\nabla \sqrt{\tilde{c}}\|_{L^3(Q_{r_k})}^3 + C r_{N+1}^{-4} \|\nabla \sqrt{\tilde{c}}\|_{L^3(Q_{r_N})}^3,
\eeno
which is similar as $I_4$ and $I_6$.
Then
\beno
I_7 \leq  CC_0^\frac32 \varepsilon_0^\frac34.
\eeno

{\bf \underline{Estimate of $I_8$}. }
For this term, noting that Proposition \ref{prop:1} and (\ref{eq:condition-fixed}), we have
\beno
I_8 \leq C\Lambda_0 \int_{Q_{r_3}} |u|^2 \leq C \Lambda_0 \varepsilon_0.
\eeno

{\bf \underline{Estimate of $I_9$}. }
Using Proposition \ref{prop:1} and (\ref{ine:induction}), we have
\beno
I_9 &\leq & \Lambda_0\sum_{k=1}^N \int_{Q_{r_k} \setminus Q_{r_{k+1}}} |u|^2\left| u \cdot \nabla \psi\right| +\Lambda_0 \int_{Q_{r_{N+1}}} |u|^2 \left|u \cdot \nabla \psi \right|\\
&\leq& C \Lambda_0\sum_{k=1}^N r_k^{-4} \|u\|_{L^\frac{10}3(Q_{r_k})}^3 r_k^\frac12 + C \Lambda_0r_{N+1}^{-4} \|u\|_{L^\frac{10}3(Q_{r_N})}^3 r_{N+1}^\frac12 \\
&\leq& C \Lambda_0\sum_{k=1}^N r_k^{-4} r_k^\frac92 C_0^\frac32 \varepsilon_0^\frac34 r_k^\frac12 + C\Lambda_0 r_{N+1}^{-4} r_N^\frac92 C_0^\frac32 \varepsilon_0^\frac34 r_{N+1}^\frac12 \\
&\leq& C \Lambda_0 C_0^\frac32 \varepsilon_0^\frac34.
\eeno

{\bf \underline{Estimate of $I_{10}$}. }
First of all, we estimate the decomposition of pressure $p $. For $0<2r<\rho\leq 1$, let $\eta\geq 0$ be supported in $B_\rho$ with $\eta\equiv1$ in $B_{\frac{\rho}2}$. The divergence of $(\ref{eq:KS})_3$ gives
$-\Delta p=\partial_i\partial_j(u_iu_j)+\nabla\cdot (n\nabla\phi)$ in the sense of distribution. Let
\ben\label{ine:p 1 u u}
p_1=\int_{\mathbb{R}^3}\frac1{4\pi|x-y|}\left[\partial_i\partial_j\left((u_i-(u_i)_\rho)(u_j-(u_j)_\rho)\eta\right]+\nabla\cdot(n\nabla\phi\eta)\right](y,t)dy
\een
where $(u)_\rho$ denotes the mean value of $u $ in $B_\rho$.
Then $p_2=p-{p_1}$ and $\Delta p_2 = 0$ in $B_\frac \rho2$. 
Choose $\rho=1$ from now on.

By Lemma \ref{mean value property} and H\"{o}lder inequality,
\ben\label{ine:harmonic estimate p 2}
\int_{{B_r}}|p_2-(p_2)_{B_r}|^{\frac32}dx &\leq& Cr^{\frac32}\int_{{B_r}}|\nabla p_2|^{\frac32}dx\\\nonumber
&\leq& C \left(\frac r\rho\right)^{\frac92} \int_{B_{\rho}}|p|^{\frac32}dx+C \left(\frac r\rho\right)^{\frac92}\int_{B_{\rho}}|p_1|^{\frac32}dx
\een
By the Calderon-Zygmund estimate and Riesz potential estimate, we have
\ben\label{ine:CZ estimate p 1}
\int_{B_{\rho}}|p_1|^{\frac32}dx\leq C \int_{B_{\rho}}|u-(u)_{B_\rho}|^3+C \rho^{\frac 34}\left(\int_{B_{\rho}}|n\nabla \phi|^\frac65dx\right)^{\frac54}.
\een

Second,
we choose $\chi_k$ to be a cut-off function which vanishes outside of
$Q_{r_k}$, equals 1 in $Q_{\frac78r_k}$ and satisfies $|\nabla\chi_k|\leq C r_k^{-1}$. Noting that $\chi_1\phi_n=\phi_n.$
\beno
I_{10}&= &C\Lambda_0\left|\int_{-1}^t\int_{B_{r_3}} p u \cdot \nabla \psi \right|\\
&\leq& C\Lambda_0 \left|\sum_{k=1}^N\int_{{Q_{r_3}}^t} (p-(p)_{B_{r_k}}) u \cdot \nabla ((\chi_k-\chi_{k+1})\psi) \right|\\
&&+C\Lambda_0 \left|\int_{{Q_{r_3}}^t} (p-(p)_{B_{r_{N+1}}}) u \cdot \nabla (\chi_{N+1}\psi) \right|\\
&\leq& C\Lambda_0(I_{10}'+I_{10}'')
\eeno
where
\beno
I_{10}'&=& \left|\sum_{k=1}^N\int_{{Q_{r_3}}^t}(p_1-(p_1)_{B_{r_k}}) u \cdot \nabla ((\chi_k-\chi_{k+1})\psi) \right|+\left| \int_{{Q_{r_3}}^t}(p_1-(p_1)_{B_{r_{N+1}}}) u \cdot \nabla (\chi_{N+1}\psi) \right|\\
\eeno
and
\beno
I_{10}''= \left| \sum_{k=1}^N\int_{{Q_{r_3}}^t}(p_2-(p_2)_{B_{r_k}}) u \cdot \nabla ((\chi_k-\chi_{k+1})\psi) \right|+ \left|\int_{{Q_{r_3}}^t}(p_2-(p_2)_{B_{r_{N+1}}}) u \cdot \nabla (\chi_{N+1}\psi) \right|.
\eeno
For $I_{10}'$, there holds
\beno
I_{10}' &\leq & \left|\sum_{k=1}^N \int_{{Q_{r_k}}^t \setminus {Q_{r_{k+2}}}^t}(p_1-(p_1)_{B_{r_k}}) u \cdot \nabla ((\chi_k-\chi_{k+1})\psi)\right| +  \left| \int_{{Q_{r_{N+1}}}^t} (p_1-(p_1)_{B_{r_{N+1}}}) u \cdot \nabla (\chi_{N+1}\psi)\right|\\
&:=& T_1'+T_2'.
\eeno
In order to estimate $T_1'$, we introduce a new cut-off function $\xi_\ell(x) = \eta(\frac{x}{r_\ell })$ and $\xi_0 = 1$. Then by $\eqref{ine:p 1 u u}$ and $1 = \sum_{\ell=0}^k (\xi_\ell - \xi_{\ell+1})+ \xi_{k+1}$, in $B_1$ for $1\leq k\leq N$, we have
\ben\label{ine:M1 M2} \nonumber
T_1'&=& \sum_{k=1}^N \int_{Q_{r_k} \setminus Q_{r_{k+2}}} \left|\int_{\mathbb{R}^3}\frac1{4\pi|x-y|}{\left(\partial_i\partial_j(u_i u_j\eta)+\nabla\cdot(n\nabla\phi\eta)\right)}(y,t)dy ~~ u \cdot \nabla ((\chi_k-\chi_{k+1})\psi)\right| \\ \nonumber
&=& \sum_{k=4}^N \int_{Q_{r_k} \setminus Q_{r_{k+2}}} \left|\int_{\mathbb{R}^3}\frac1{4\pi|x-y|}{\partial_i\partial_j\left(u_iu_j\left[\sum_{\ell=0}^{k} (\xi_\ell - \xi_{\ell+1})+ \xi_{k+1}\right]\eta\right)}(y,t)dy u \cdot \nabla ((\chi_k-\chi_{k+1})\psi)\right|\\ \nonumber
&+&\sum_{k=4}^N \int_{Q_{r_k} \setminus Q_{r_{k+2}}} \left|\int_{\mathbb{R}^3}\frac1{4\pi|x-y|}{\nabla\cdot\left(n\nabla\phi\left[\sum_{\ell=0}^{k} (\xi_\ell - \xi_{\ell+1})+ \xi_{k+1}\right]\eta\right)}(y,t)dy u \cdot \nabla ((\chi_k-\chi_{k+1})\psi)\right| \\ \nonumber&&+
\sum_{k=1}^3 \int_{Q_{r_k} \setminus Q_{r_{k+2}}} \left|\int_{\mathbb{R}^3}\frac1{4\pi|x-y|}{\left(\partial_i\partial_j(u_i u_j\eta)+\nabla\cdot(n\nabla\phi\eta)\right)}(y,t)dy ~~ u \cdot \nabla ((\chi_k-\chi_{k+1})\psi)\right|\\
&:=& M_1 + M_2+M_3,
\een
Noting that the support set of $\xi_\ell$, we have
\beno
M_1 &=& \sum_{k=4}^N \int_{Q_{r_k} \setminus Q_{r_{k+2}}} \left|\sum_{\ell=0}^{k-3} \int_{B_{r_\ell} \setminus B_{r_{\ell+2}}}\frac1{4\pi|x-y|}{\left(\partial_i\partial_j\left(u_iu_j (\xi_\ell - \xi_{\ell+1})\eta\right)\right)}(y,t)dy u \cdot \nabla ((\chi_k-\chi_{k+1})\psi)\right| \\
&&+ \sum_{k=4}^N \int_{Q_{r_k} \setminus Q_{r_{k+2}}} \left|\int_{B_{r_{k-2}}}\frac1{4\pi|x-y|}{\left(\partial_i\partial_j\left(u_iu_j \xi_{k-2}\eta\right)\right)}(y,t)dy u \cdot \nabla ((\chi_k-\chi_{k+1})\psi)\right| \\
&:=& M_{11} + M_{12} .
\eeno
Since $|x-y| \geq r_{\ell+3}$ for the term $M_{11}$ and $|\nabla ((\chi_k-\chi_{k+1})\psi)| \leq C r_k^{-4}$, by using integration by parts, H\"{o}lder's inequality and $\eqref{ine:induction}$, there holds
\ben\label{ine:M11} \nonumber
M_{11} &\leq& C \sum_{k=4}^N \int_{Q_{r_k} \setminus Q_{r_{k+2}}} \left|\sum_{\ell=0}^{k-3} \int_{B_{r_\ell}}r_{\ell+3}^{-3}{\left(u_iu_j (\xi_\ell - \xi_{\ell+1})\eta\right)}(y,t)dy u \cdot \nabla ((\chi_k-\chi_{k+1})\psi)\right| \\ \nonumber
&\leq& C \sum_{k=4}^N r_k^{-4} \int_{Q_{r_k}} \left|\sum_{\ell=0}^{k-3} \int_{B_{r_\ell}}r_{\ell+3}^{-3}|u|^2dy\right| |u|  \\ \nonumber
&\leq& C \sum_{k=4}^N r_k^{-4} \|u\|_{L^3(Q_{r_k})} \left(\int_{Q_{r_k}} \left|\sum_{\ell=0}^{k-3} \int_{B_{r_\ell}}r_{\ell+3}^{-3}|u|^2dy\right|^\frac32\right)^\frac23 \\ \nonumber
&\leq& C \sum_{k=4}^N r_k^{-4} r_k^\frac53 C_0^\frac12 \varepsilon_0^\frac14 r_k^\frac{10}3 \left|\sum_{\ell=0}^{k-3} r_{\ell+3}^{-3} r_k^2r_{\ell} C_0 \varepsilon_0^\frac12\right| \\\nonumber
&\leq& C \sum_{k=4}^N r_k^{-4} r_k^\frac53 C_0^\frac12 \varepsilon_0^\frac14 r_k^\frac{10}3 \left|\sum_{\ell=0}^{k-3} r_{\ell+3}^{-3} r_{\ell}^3 C_0 \varepsilon_0^\frac12\right| \\
&\leq& C \sum_{k=4}^N  k r_k C_0^\frac32 \varepsilon_0^\frac34 \leq C C_0^\frac32 \varepsilon_0^\frac34.
\een
For the term $M_{12}$, noting that $|x-y|=0$ for some $y \in B_{r_k}$ and $(x,t) \in Q_{r_k} \setminus Q_{r_{k+2}}$, the method of the estimate of $M_{11}$ is fail. Noting that the operator $T$ operates on any function $F$, and $TF=\int_{\mathbb{R}^3} |x-y|^{-1} \partial_i \partial_j F_{ij}$, $T$ satisfies the conditions of 
singular integral theorem, and by H\"{o}lder's inequality, there holds
\ben\label{ine:M12} \nonumber
M_{12} &\leq& C \sum_{k=4}^N r_k^{-4} \|u\|_{L^3(Q_{r_k})} \left(\int_{Q_{r_k}} \left|\int_{B_{r_{k-2}}}\frac1{4\pi|x-y|}{\left(\partial_i\partial_j\left(u_iu_j \xi_{k-2}\eta\right)\right)}(y,t)dy\right|^\frac32\right)^\frac23 \\ \nonumber
&\leq& C \sum_{k=4}^N r_k^{-4} \|u\|_{L^3(Q_{r_k})}\left(\int_{-r_k^2}^0 \||u|^2 \xi_{k-2}\eta\|^{\frac32}_{L^\frac32(\mathbb{R}^3)} dt\right)^{\frac23} \\
&\leq& C \sum_{k=4}^N r_k^{-4} \|u\|_{L^3(Q_{r_k})} \|u\|_{L^3(Q_{r_{k}})}^2 \leq C C_0^\frac32 \varepsilon_0^\frac34.
\een
Collecting $\eqref{ine:M11}$, $\eqref{ine:M12}$ ,and we have
\ben\label{ine:M1}
M_1 \leq C C_0^\frac32 \varepsilon_0^\frac34.
\een

The estimate of the term $M_2$ is same as $M_1$,
\beno
M_2 &=& \sum_{k=4}^N \int_{Q_{r_k} \setminus Q_{r_{k+2}}} \left|\sum_{\ell=0}^{k-3} \int_{B_{r_\ell} \setminus B_{r_{\ell+2}}}\frac1{4\pi|x-y|}\nabla \cdot (n \nabla \phi (\xi_\ell-\xi_{\ell+1}) \eta)(y,t)dy u \cdot \nabla ((\chi_k-\chi_{k+1})\psi)\right| \\
&&+ \sum_{k=4}^N \int_{Q_{r_k} \setminus Q_{r_{k+2}}} \left|\int_{B_{r_{k-2}}}\frac1{4\pi|x-y|}{\nabla \cdot (n \nabla \phi \xi_{k-2} \eta)}(y,t)dy u \cdot \nabla ((\chi_k-\chi_{k+1})\psi)\right| \\
&:=& M_{21} + M_{22} .
\eeno
For $M_{21}$, since $|x-y| \geq r_{\ell+3}$, using $\eqref{ine:induction}$ and H\"{o}lder's inequality, there holds
\ben\label{ine:M21} \nonumber
M_{21} &\leq& \sum_{k=4}^N \int_{Q_{r_k} \setminus Q_{r_{k+2}}} \left|\sum_{\ell=0}^{k-3} \int_{B_{r_\ell} \setminus B_{r_{\ell+2}}} r_{\ell+3}^{-2} |\nabla \phi| |n|dy u \cdot \nabla ((\chi_k-\chi_{k+1})\psi)\right| \\ \nonumber
&\leq& C \sum_{k=4}^N r_k^{-4} \int_{Q_{r_k}} \left|\sum_{\ell=0}^{k-3}\int_{B_{r_\ell}} r_{\ell+3}^{-2} |\nabla \phi| |n|dy |u| \right|\\ \nonumber
&\leq& C \sum_{k=4}^N r_k^{-4} \int_{Q_{r_k}} \sum_{\ell=0}^{k-3} r_{\ell+3}^{-2} r_{\ell}^3 C_0 \varepsilon_0^\frac12 |u| \\
&\leq& C \sum_{k=4}^N r_k^{-4} r_k^\frac{10}3 \|u\|_{L^3(Q_{r_k})} C_0 \varepsilon_0^\frac12 \leq C C_0^\frac32 \varepsilon_0^\frac34.
\een
For the term $M_{22}$, by H\"{o}lder's inequality, Riesz potential estimate and the condition $\eqref{ine:induction}$, there holds
\ben\label{ine:M22} \nonumber
M_{22} &\leq& \sum_{k=4}^N \int_{Q_{r_k} \setminus Q_{r_{k+2}}} \left|\int_{B_{k-2}}\frac1{4\pi|x-y|} \nabla \cdot (n \nabla \phi \xi_{k-2} \eta)(y,t) dy u \cdot \nabla ((\chi_k-\chi_{k+1})\psi)\right| \\ \nonumber
&\leq& C \sum_{k=4}^N \int_{Q_{r_k} \setminus Q_{r_{k+2}}} \left|\int_{B_{r_{k-2}}} |x-y|^{-2} |n \nabla \phi \xi_{k-2} \eta|(y,t) dy u \cdot \nabla ((\chi_k-\chi_{k+1})\psi)\right| \\ \nonumber
&\leq& C \sum_{k=4}^N r_k^{-4} \|u\|_{L^3(Q_{r_k})} \left(\int_{-r_k^2}^0 \left(\int_{B_{r_k}} \left|\int_{B_{r_{k-2}}} |x-y|^{-2} |n \nabla \phi \xi_{k-2} \eta|(y,t) dy\right|^\frac{15}4dx \right)^\frac4{9}dt\right)^\frac35 r_k^\frac43 \\ \nonumber
&\leq& C \sum_{k=4}^N r_k^{-\frac83} \|u\|_{L^3(Q_{r_k})} \left(\int_{-r_k^2}^0 \||n \nabla \phi \xi_{k-2} \eta|\|_{L^\frac53(\mathbb{R}^3)}^\frac53 \right)^\frac35 \\ \nonumber
&\leq& C \sum_{k=4}^N r_k^{-\frac83} \|u\|_{L^3(Q_{r_k})} \|n\|_{L^\frac53(Q_{r_k})} \\
&\leq& C \sum_{k=4}^N r_k^{-\frac83} r_k^\frac53 C_0^\frac12 \varepsilon_0^\frac14 r_k^3 C_0 \varepsilon_0^\frac12 \leq C C_0^\frac32 \varepsilon_0^\frac34.
\een
Collecting $\eqref{ine:M21}$, $\eqref{ine:M22}$, we have
\ben\label{ine:M2}
M_2 \leq C C_0^\frac32 \varepsilon_0^\frac34.
\een
As for the term $M_3$,  by H\"{o}lder inequality there holds
\ben\label{ine:M3}
&M_3 &\nonumber
=\sum_{k=1}^3 \int_{Q_{r_k} \setminus Q_{r_{k+2}}} p_1 u \cdot \nabla ((\chi_k-\chi_{k+1})\psi)\\
&&\leq C\sum_{k=1}^3r_k^{-4}\|u\|_{L^3(Q_{r_k})}\|p_1\|_{L^{\frac32}(Q_{r_k})} \leq C \|u\|_{L^3(Q_{r_1})}\|p_1\|_{L^{\frac32}(Q_{r_1})}.
\een
By (\ref{eq:condition-fixed}), (\ref{ine:CZ estimate p 1}) and (\ref{ine:M3}), we arrive
\ben\label{ine:M,3,} \nonumber
M_3 &\leq& C \|u\|_{L^3(Q_{r_1})}\left(C\int_{I_{r_1}} \left[\int_{B_{r_1}}|u-(u)_{B_{r_1}}|^3+ \left(\int_{B_{r_1}}|n\nabla \phi|^\frac65dx\right)^{\frac54} \right]dt\right)^{\frac23}\\\nonumber
&\leq& C \|u\|^3_{L^3(Q_{r_1})} + C \|\nabla \phi\|_{L^\infty} \|u\|_{L^3(Q_{r_1})} \|n\|_{L^\frac53(Q_{r_1})}^\frac23 \\
&\leq& (\|\nabla  \phi\|_{L^\infty}+1) C (\varepsilon_0^\frac32 + \varepsilon_0^\frac56).
\een

To sum up, $\eqref{ine:M1 M2}$, $\eqref{ine:M1}$, $(\ref{ine:M2})$  and $\eqref{ine:M,3,}$ implies that
\beno
T_1' \leq C C_0^\frac32 \varepsilon_0^\frac34+ C \varepsilon_0^\frac56 \leq C C_0^\frac32 \varepsilon_0^\frac34.
\eeno
The term of $T_2'$ is same as $T_1'$, we omit the estimate of $T_2'$ and arrive
\beno
I_{10}' \leq C C_0^\frac32 \varepsilon_0^\frac34.
\eeno
For the another term $I_{10}''$, by (\ref {mean}), we have 
\beno
I_{10}'' &\leq & \sum_{k=1}^N \int_{Q_{r_k} \setminus Q_{r_{k+2}}}\left| (p_2 - (p_2)_{B_{r_k}}) u \cdot  \nabla ((\chi_k-\chi_{k+1})\psi)\right|\\
 &&+ \int_{Q_{r_{N+1}}} \left|(p_2 - (p_2)_{B_{r_{N}}}) u \cdot \nabla (\chi_{N+1}\psi)\right|\\
&\leq& C \sum_{k=1}^N r_k^{-4} \|p_2 - (p_2)_{B_{r_k}}\|_{L^\frac32(Q_{r_k})} \|u\|_{L^3(Q_{r_k})} \\
&&+ C  r_{N+1}^{-4} \|p_2 - (p_2)_{B_{N}}\|_{L^\frac32(Q_{r_{N}})} \|u\|_{L^3(Q_{r_{N}})}.
\eeno
Using (\ref{ine:CZ estimate p 1}), (\ref{ine:harmonic estimate p 2}) and (\ref{eq:condition-fixed}), for any $k = 1,2,\cdots,$ we have
\beno
&&\|p_2-(p_2)_{B_{r_k}}\|_{L^\frac32(Q_{r_k})}^\frac32\\ &\leq& C r_k^\frac92 \|p_2\|_{L^\frac32(Q_{r_1})}^\frac32 \\
&\leq& C r_k^\frac92 \|p_1\|_{L^\frac32(Q_{r_1})}^\frac32 + C r_k^\frac92 \|p\|_{L^\frac32(Q_{r_1})}^\frac32 \\
&\leq& C r_k^\frac92 \|u\|_{L^3(Q_{r_1})}^3 + C r_k^\frac92 \int_{I_{r_1}} \left(\int_{B_{r_1}}|n\nabla \phi|^\frac65dx\right)^{\frac54} dt + C r_k^\frac92 \|p\|_{L^\frac32(Q_{r_1})}^\frac32 \\
&\leq& C r_k^\frac92 (\|\nabla\phi\|_\infty+1)^{\frac32}\varepsilon_0.
\eeno
Then
\beno
I_{10}'' \leq C \Lambda_1\sum_{k=1}^N r_k^{-4} r_k^3 \varepsilon_0^{\frac23} r_k^\frac53 C_0^\frac12 \varepsilon_0^\frac14 + C \Lambda_1r_{N+1}^{-4} r_N^3 \varepsilon_0^{\frac23} r_N^\frac53 C_0^\frac12 \varepsilon_0^\frac14 \leq C C_0^\frac12 \Lambda_1\varepsilon_0^\frac{11}{12}.
\eeno
Collecting $I_{10}'$ and $I_{10}''$, we have
\beno
I_{10} \leq C\Lambda_0 \Lambda_1 C_0^\frac32 \varepsilon_0^\frac34 + C \Lambda_0 \Lambda_1C_0^\frac12 \varepsilon_0^\frac{11}{12}\leq C  \Lambda_0\Lambda_1C_0^\frac32 \varepsilon_0^\frac34.
\eeno

{\bf \underline{Estimate of $I_{11}$}. }
Noting that $\nabla \phi \in L^\infty(Q_{r_1})$ and using Proposition \ref{prop:1}, (\ref{ine:n 5 3}) and (\ref{ine:induction}), we have
\beno
I_{11} &\leq &C \Lambda_0 \sum_{k=1}^N \int_{Q_{r_k} \setminus Q_{r_{k+1}}}\left| n \nabla \phi \cdot u \psi\right| +C \Lambda_0 \int_{Q_{r_{N+1}}}\left| n \nabla \phi \cdot u \psi \right|\\
&\leq& C \Lambda_0\Lambda_1\sum_{k=1}^N r_k^{-3} \|n\|_{L^\frac53(Q_{r_k})} \|u\|_{L^\frac{10}3(Q_{r_k})} r_k^\frac12 + C\Lambda_0\Lambda_1 r_{N+1}^{-3} \|n\|_{L^\frac53(Q_{r_N})} \|u\|_{L^\frac{10}3(Q_{r_N})} r_N^{\frac12} \\
&\leq& C\Lambda_0\Lambda_1 \sum_{k=1}^N r_k^{-3} r_k^3 C_0 \varepsilon_0^\frac12 r_k^\frac32 C_0^\frac12 \varepsilon_0^\frac14 + C \Lambda_0\Lambda_1 r_{N+1}^{-3} r_N^3 C_0 \varepsilon_0^\frac12 r_N^\frac32 C_0^\frac12 \varepsilon_0^\frac14 \\
&\leq& C \Lambda_0\Lambda_1C_0^\frac32 \varepsilon_0^\frac34.
\eeno

Collecting $I_1$, $I_2$, $\cdots$, $I_{11}$, for any $t \in (-r_{N+1}^2,0)$, we arrive at
\ben\label{ine:induction 1} \nonumber
&&\int_{B_{r_{N+1}}} (n \ln n \psi)(\cdot,t) + r_{N+1}^{-3}\int_{{Q_{r_{N+1}}}^t} |\nabla \sqrt{n}|^2 + r_{N+1}^{-3}\int_{B_{r_{N+1}}} |\nabla \sqrt{\tilde{c}}|^2\\
&& + r_{N+1}^{-3}\int_{{Q_{r_{N+1}}}^t} |\nabla^2 \sqrt{\tilde{c}}|^2+ \Lambda_0  r_{N+1}^{-3}\int_{B_{r_{N+1}}} |u|^2 + \Lambda_0  r_{N+1}^{-3}\int_{{Q_{r_{N+1}}}^t} |\nabla u|^2 \nonumber\\
&\leq& C\Lambda_0\Lambda_1 C_0^\frac32 \varepsilon_0^\frac34 +  C C_0^\frac{29}{18} \varepsilon_0^\frac{25}{36},
\een
where $\Lambda_0\geq 1.$

\textcolor[rgb]{0.00,0.00,0.00}{{\bf Step II. Estimate of  the term $r_{N+1}^{-3} \int_{B_{r_{N+1}}} n$. }}

Recall the equation of $n$ as follows:
\beno
\partial_t n - \Delta n + u\cdot \nabla n = - \nabla \cdot (n \nabla c).
\eeno
Multiplying it with $\psi$ and integration by parts on ${Q_1}^t$, we arrive at
%

%
%
%
%
%
\ben\label{ine:energy n <x>} \nonumber
\int_{B_1} (n \psi)(\cdot,t) = \int_{{Q_1}^t} n  (\partial_t \psi + \Delta \psi)
+ \int_{{Q_1}^t} n  u \cdot \nabla \psi+ \int_{{Q_1}^t} n  \nabla c \cdot \nabla \psi.
\een
Recall $\psi = \phi_{N+1}$ and using Proposition \ref{prop:1}, we arrive at

\beno
\int_{B_{r_{N+1}}} (n \psi)(\cdot,t) = \int_{{Q_{r_3}}^t} n  (\partial_t \psi + \Delta \psi)
+ \int_{{Q_{r_3}}^t} n  u \cdot \nabla \psi+ \int_{{Q_{r_3}}^t} n  \nabla c \cdot \nabla \psi=: K_1+  \cdots + K_3.
\eeno

Using (\ref{eq:condition-fixed}) and Proposition \ref{prop:1}, for $K_1$, we have
\beno
K_1 \leq C \int_{{Q_{r_3}}^t} n  \leq C \varepsilon_0.
\eeno

For $K_2$, using Proposition \ref{prop:1}, direct calculations imply that
\beno
K_2 &\leq & \sum_{k=1}^N \int_{Q_{r_k} \setminus Q_{r_{k+1}}} nu\cdot \nabla \psi + \int_{Q_{r_{N+1}}}nu\cdot \nabla \psi \\
&\leq& C \sum_{k=1}^N r_k^{-4}\|n\|_{L^\frac53(Q_{r_k})} \|u\|_{L^\frac{10}3(Q_{r_k})} r_k^\frac12 + C \|n\|_{L^\frac53(Q_{r_N})}r_N^{-4} \|u\|_{L^\frac{10}3(Q_{r_N})} r_N^\frac12 \\
&\leq& C \sum_{k=1}^N r_k^{-4}r_k^3C_0\varepsilon_0^\frac12 r_k^\frac32 C_0^\frac12 \varepsilon_0^\frac14 r_k^\frac12+Cr_N^{-4}r_N^3C_0\varepsilon_0^\frac12 r_N^\frac32 C_0^\frac12 \varepsilon_0^\frac14 r_N^\frac12\\
&\leq& C C_0^\frac32 \varepsilon_0^\frac34.
\eeno

The term $K_3$ is similar as the term $K_2$. Noting that $|\nabla c|\leq 2|\sqrt{\tilde{c}} \nabla \sqrt{\tilde{c}}|\leq 2\sqrt{\Lambda_0}|\nabla \sqrt{\tilde{c}}|$, using Proposition \ref{prop:1}, we have
\beno
K_3 \leq C \sqrt{\Lambda_0} C_0^\frac32 \varepsilon_0^\frac34.
\eeno

Collecting the estimates of $K_1-K_3$, for any $t \in (-r_{N+1}^2,0)$, we have
\ben\label{ine:induction 2}
\int_{B_{r_{N+1}}} (n\psi)(\cdot,t) \leq C\varepsilon_0 + C C_0^\frac32 \varepsilon_0^\frac34+C \sqrt{\Lambda_0} C_0^\frac32 \varepsilon_0^\frac34 \leq   C\sqrt{\Lambda_0} C_0^\frac32 \varepsilon_0^\frac34.
\een

{\bf Step III. Estimate of $r_{N+1}^{-3} \int_{B_{r_{N+1}}} n|\ln n|$.}\\

From step II, we know that in order to prove (\ref{ine:induction}) for $k = N+1$, it's sufficient to estimate the term of  $r_{N+1}^{-3} \int_{B_{r_{N+1}}} n + n |\ln n|$. Let's prove it. First, combining (\ref{ine:induction 1}) and (\ref{ine:induction 2}), for any $t \in (-r_{N+1}^2,0)$, we have
\ben\label{eq:n ln estimate} \nonumber
&&\int_{B_{r_{N+1}}} (n \psi)(\cdot,t)+\int_{B_{r_{N+1}}} (n \ln n \psi)(\cdot,t) + r_{N+1}^{-3}\int_{{Q_{r_{N+1}}}^t} |\nabla \sqrt{n}|^2 \\
&&+ r_{N+1}^{-3}\int_{B_{r_{N+1}}} (|\nabla \sqrt{\tilde{c}}|^2)(\cdot,t) + r_{N+1}^{-3}\int_{{Q_{r_{N+1}}}^t} |\nabla^2 \sqrt{\tilde{c}}|^2 \nonumber \\
&&+ \Lambda_0  r_{N+1}^{-3}\int_{B_{r_{N+1}}} (|u|^2)(\cdot,t) + \Lambda_0  r_{N+1}^{-3}\int_{{Q_{r_{N+1}}}^t} |\nabla u|^2 \nonumber\\
&\leq& C C_0^{\frac{9}{10}} \varepsilon_0^{\frac{11}{20}}+C\Lambda_0 \Lambda_1C_0^\frac32 \varepsilon_0^\frac34  +  C C_0^\frac{29}{18} \varepsilon_0^\frac{25}{36}\leq  C C_0^{\frac{29}{18}} \varepsilon_0^{\frac{11}{20}}+C\Lambda_0\Lambda_1 C_0^\frac32 \varepsilon_0^\frac34
\een



Using $|\ln n| n^{\alpha}\leq \alpha^{-1}e^{-1}$ for $0<n<1$ and $0<\alpha<\frac{1}{20}$, by (\ref{ine:induction 2}) and the above inequality, we have
\ben\label{nlnn'} \nonumber
&&Cr_{N+1}^{-3} \int_{B_{r_{N+1}}} (n |\ln n|)(\cdot,t) dx\\
&\leq& \int_{B_{r_{N+1}}} (n\ln n \psi)(\cdot,t) dx -2\int_{B_{r_{N+1}}\cap \{x;0<n<1\}}(n\ln n \psi)(\cdot,t) dx
\nonumber \\
&\leq &\int_{B_{r_{N+1}}} (n\ln n \psi)(\cdot,t) dx + 2\alpha^{-1}e^{-1} \int_{B_{r_{N+1}}}  (n^{1-\alpha}\psi)(\cdot,t) dx
\nonumber \\\nonumber
&\leq& C C_0^{\frac{29}{18}} \varepsilon_0^{\frac{11}{20}}+C\Lambda_0\Lambda_1 C_0^\frac32 \varepsilon_0^\frac34 +2\alpha^{-1}e^{-1}\left(C C_0^{\frac{9}{10}} \varepsilon_0^{\frac{11}{20}} +   C\sqrt{\Lambda_0} C_0^\frac32 \varepsilon_0^\frac34\right)^{1-\alpha}\\
&\leq&2\alpha^{-1}e^{-1}\left(C C_0^{\frac{29}{18}} \varepsilon_0^{\frac{11}{20}}+C\Lambda_0 \Lambda_1C_0^\frac32 \varepsilon_0^\frac34\right)^{1-\alpha}
\een
where we used the integral of heat kernel is bounded
\beno
\int_{B_{r_{N+1}}}  \psi dx\leq C,
\eeno
and we will choose the smallness of the right term $C C_0^{\frac{29}{18}} \varepsilon_0^{\frac{11}{20}}+C\Lambda_0 \Lambda_1C_0^\frac32 \varepsilon_0^\frac34$.

%
%
%
{\bf Step IV. The proof of (\ref{ine:induction}). }
Combining (\ref{ine:induction 2})-(\ref{nlnn'}), for any $t \in (-r_{N+1}^2,0)$, we arrive at
\beno
&&r_{N+1}^{-3} \int_{B_{r_{N+1}}} (n + |n \ln n| + |\nabla \sqrt{\tilde{c}}|^2 + |u|^2)(\cdot,t) \\
&&+ r_{N+1}^{-3} \int_{{Q_{r_{N+1}}}^t} |\nabla \sqrt{n}|^2 + |\nabla u|^2 + |\nabla^2 \sqrt{\tilde{c}}|^2 \\
&\leq&\alpha^{-1}\left(C C_0^{\frac{29}{18}} \varepsilon_0^{\frac{11}{20}}+C\Lambda_0 \Lambda_1C_0^\frac32 \varepsilon_0^\frac34\right)^{1-\alpha}\\
&\leq & \left(C\alpha^{-1}C_0^2\varepsilon_0^{\frac{1}{20}-\frac{11}{20}\alpha}+ C\alpha^{-1}C_0^2(\Lambda_0 \Lambda_1)^{1-\alpha}\varepsilon_0^{\frac{1}{4}-\frac{3}{4}\alpha}\right)C_0\varepsilon_0^{\frac{1}{2}}
\eeno
Due to (\ref{eq:condition-fixed}) and $\varepsilon_0=\frac{\varepsilon_1}{(\Lambda_0 \Lambda_1)^{4+4\alpha_0}}$, we choose  $\alpha=\min\{\frac{1}{20},\frac{\alpha_0}{4+6\alpha_0}\}$, and $\varepsilon_1$ such that
\beno
C\alpha^{-1}C_0^2\varepsilon_0^{\frac{1}{20}-\frac{11}{20}\alpha}\leq C\alpha^{-1}C_0^2\varepsilon_1^{\frac{1}{20}-\frac{11}{20}\alpha}\leq \frac12
\eeno
and
\beno
C\alpha^{-1}C_0^2(\Lambda_0 \Lambda_1)^{1-\alpha-(\frac{1}{4}-\frac{3}{4}\alpha)(4+4\alpha_0)}\varepsilon_1^{\frac{1}{4}-\frac{3}{4}\alpha}\leq \frac12
\eeno
since
\beno
1-\alpha-(\frac{1}{4}-\frac{3}{4}\alpha)(4+4\alpha_0)<0\Leftrightarrow 0<\alpha<\frac{\alpha_0}{2+3\alpha_0}.
\eeno

Then we have
\beno
&&r_{N+1}^{-3} \int_{B_{r_{N+1}}} (n+ |n \ln n| + |\nabla \sqrt{\tilde{c}}|^2 + |u|^2)(\cdot,t) \\
&&\quad+ r_{N+1}^{-3} \int_{{Q_{r_{N+1}}}^t} |\nabla \sqrt{n}|^2 + |\nabla u|^2 + |\nabla^2 \sqrt{\tilde{c}}|^2  \leq C_0 \varepsilon_0^\frac12.
\eeno

Then, for uniform $t \in (-r_{N+1}^2,0)$, we have
\beno
&&r_{N+1}^{-3} \sup_{-r_{N+1}^2<t<0} \int_{B_{r_{N+1}}} n+ |n \ln n| + |\nabla \sqrt{\tilde{c}}|^2 + |u|^2 \\
&&\quad+ r_{N+1}^{-3} \int_{Q_{r_{N+1}}} |\nabla \sqrt{n}|^2 + |\nabla u|^2 + |\nabla^2 \sqrt{\tilde{c}}|^2  \leq C_0 \varepsilon_0^\frac12.
\eeno

To sum up, for any $k=1,2,\cdots,N+1,N+2,\cdots$, (\ref{ine:induction}) is true.

{\bf Step V. The proof of bounded-ness of $(n,\nabla c,u)$.}

In the end, by interpolation inequality, we achieve at
\beno
\left(\int_{Q_{r_k}}|\sqrt{n}|^{\frac{10}3}dxdt\right)^{\frac3{10}}\leq C || \sqrt{n}||^{\frac25}_{L^\infty_tL^2_x(Q_{r_k})} ||\nabla \sqrt{n}||^{\frac35}_{L^2_tL^2_x(Q_{r_k})}+ C || \sqrt{n}||_{L^\infty_tL^2_x(Q_{r_k})}.
\eeno
and
\beno
\left(r_k^{-5}\int_{Q_{r_k}}|\sqrt{n}|^{\frac{10}3}dxdt\right)^{\frac3{10}}\leq C_0^{\frac15+\frac3{10}}r_k^{-\frac32}  r_k^{\frac35}\varepsilon_0^{\frac1{10}}r_k^\frac{9}{10}\varepsilon_0^{\frac3{20}}\leq C_0^\frac12\varepsilon_0^{\frac14}.
\eeno
Applying Lebesgue differential theorem, there holds

%

\beno
n(0,0)\leq C_3,
\eeno
where $C_3 = C_0 \varepsilon_0^\frac12$.
%
Similarly, we have
\beno
\left(r_k^{-5}\int_{Q_{r_k}}|\nabla \sqrt{\tilde{c}}|^{\frac{10}3}dxdt\right)^{\frac3{10}}\leq C_3^\frac12,
\eeno
and
\beno
\left(r_k^{-5}\int_{Q_{r_k}}|u|^{\frac{10}3}dxdt\right)^{\frac3{10}}\leq C_3^\frac12,
\eeno
which means
\ben\label{boundeness}
|\nabla \sqrt{\tilde{c}}(0,0)|\leq C_1^\frac12; \quad |u(0,0)|\leq C_3^\frac12.
\een
Using $(\ref{boundeness})$, we have
\beno
|\nabla c (0,0)| = |\nabla \tilde{c} (0,0)|= |2\sqrt{\tilde{c}} \nabla \sqrt{\tilde{c}}| \leq C C_3^\frac12.
\eeno
Therefore, the proof is complete.

\section{The singular set's estimate}

{\bf Proof of Theorem {\ref{thm:lin}}.} It suffices to prove the inequality
(\ref{eq:condition-32}), since one can use the embedding inequality
\beno\label{ine:estimate n ln n-}
\int_{Q_{1}} |n \ln n|^{\frac32} \leq C \int_{Q_{1} \cap \{ n(x) \leq 100\}} |n|^{\frac43} + C \int_{Q_{1} \cap \{ n(x) > 100\}} |n|^\frac53.
\eeno

{\bf Step I.}
It follows from the local energy inequality that
\ben\label{eq:local energy inequality-}
&& \int_{B_1} (n \ln n \zeta)(\cdot,t) + 2 \int_{{Q_1}^t} |\nabla \sqrt{n}|^2 \zeta + \frac12   \int_{B_1} (|\nabla \sqrt{\tilde{c}}|^2 \zeta)(\cdot,t) +  \int_{{Q_1}^t} |\nabla^2 \sqrt{\tilde{c}}|^2 \zeta \nonumber\\ \nonumber
&& \quad + \frac1{36}\int_{{Q_1}^t} (\sqrt{\tilde{c}})^{-2} (\partial_j \sqrt{\tilde{c}})^2 (\partial_i \sqrt{\tilde{c}})^2 \zeta + \frac12 \int_{{Q_1}^t} |\nabla \sqrt{\tilde{c}}|^2 n \zeta \\ \nonumber
&&+\Lambda_0\int_{B_1} (|u|^2 \zeta)(\cdot,t) + \Lambda_0\int_{{Q_1}^t} |\nabla u|^2 \zeta\\ \nonumber
&\leq& \int_{{Q_1}^t} n \ln n (\partial_t \zeta + \Delta \zeta) + \int_{{Q_1}^t} n \ln n u \cdot \nabla \zeta + \int_{{Q_1}^t} n \ln n \nabla c \cdot \nabla \zeta + \int_{{Q_1}^t} n \nabla c \cdot \nabla \zeta \\ \nonumber
&&+ 2\int_{{Q_1}^t} |\nabla \sqrt{\tilde{c}}|^2 (\partial_t \zeta + \Delta \zeta) + 2 \int_{{Q_1}^t} |\nabla \sqrt{\tilde{c}}|^2 u \cdot \nabla \zeta\\
&&+\frac{20}3 \int_{{Q_1}^t} (\sqrt{\tilde{c}})^{-1} |\nabla \sqrt{\tilde{c}}|^2 |\nabla \sqrt{\tilde{c}} \cdot \nabla \zeta|
+\Lambda_0\int_{{Q_1}^t} |u|^2 (\partial_t \zeta + \Delta \zeta)+ \Lambda_0\int_{{Q_1}^t} |u|^2  u \cdot \nabla \zeta
 \nonumber\\
&&+ \Lambda_0\int_{{Q_1}^t} (p - \bar{p}) u \cdot \nabla \zeta - 2\Lambda_0\int_{{Q_1}^t} n \nabla \phi \cdot u \zeta
\een
where $t \in (-1,0)$ and $\zeta$ is a cut-off function on domain $Q_{1}$, which means that $\zeta = 1$ on $Q_{\frac12}$ and $\zeta = 0$ outside $Q_{1}$.
Moreover,
\beno
|\nabla \zeta|+|\nabla^2\zeta|+|\partial_t\zeta|\leq C.
\eeno

We just compute the first two terms, and other terms are similar.
Using H\"{o}lder inequality,
\beno\label{ine:estimate n ln n-}
A_1\leq C\int_{{Q_{1}}^t} |n \ln n|dxdt \leq C \int_{Q_1} |n \ln n| \leq C \left(\int_{Q_1} |n \ln n|^\frac32\right)^\frac23\leq C\varepsilon_2^{\frac23}.
\eeno
and
\beno
A_{2}\leq C\int_{{Q_1}^t} |n \ln n| |u|dxdt \leq C \|n \ln n\|_{L^\frac32(Q_1)} \|u\|_{L^3(Q_1)} \leq C\varepsilon_2
\eeno

Hence, we have
\ben\label{eq:nlnn}
&&\int_{B_1} (n \ln n \zeta)(\cdot,t) + \int_{{Q_1}^t} |\nabla \sqrt{n}|^2 \zeta + \int_{B_1} (|\nabla \sqrt{\tilde{c}}|^2 \zeta)(\cdot,t)\nonumber\\ \nonumber
&& \quad +  \int_{{Q_1}^t} |\nabla^2 \sqrt{\tilde{c}}|^2 \zeta+\Lambda_0\int_{B_1} (|u|^2 \zeta)(\cdot,t) + \Lambda_0\int_{{Q_1}^t} |\nabla u|^2 \zeta\\ 
&\leq& C\Lambda_0\Lambda_1\varepsilon_2^{\frac23}.
\een

{\bf Step II.} Similar as (\ref{ine:energy n <x>}), for any $t \in (-1,0)$, we get
\ben\label{ine:energy n <x>'} \nonumber
\int_{B_1} (n  \zeta)(\cdot,t) &=& \int_{{Q_1}^t} n (\partial_t \zeta + \Delta \zeta) + \int_{{Q_1}^t} n u \cdot \nabla \zeta + \int_{{Q_1}^t} n \nabla c \cdot \nabla \zeta
\een
Thus we obtain that
\ben\label{eq:n b1}
\int_{B_1} (n \zeta)(\cdot,t) \leq C\int_{{Q_1}}[n+|u|^3+n^{\frac32}+\Lambda_0 n |\nabla \sqrt{c+1}|] dxdt\leq C\Lambda_0\varepsilon_2^{\frac23}
\een

{\bf Step III.}
Recalling the estimate of  (\ref{nlnn'}), for any $t\in(-1,0)$, we have
\ben\label{nlnn''}
 \int_{B_{1}} (n |\log n|\zeta)(\cdot,t) dx&=& \int_{B_{1}} (n\log n \zeta)(\cdot,t) dx -2\int_{B_{1}\cap \{x;0<n<1\}}(n\log n \zeta)(\cdot,t) dx
\nonumber \\
&\leq &\int_{B_{1}} (n\log n \zeta)(\cdot,t) dx +C \int_{B_{1}}  (\sqrt{n}\zeta)(\cdot,t) dx.
\een
Combining this with the estimates of (\ref{eq:nlnn}), (\ref{eq:n b1}) and (\ref{nlnn''}), we have
\ben\label{eq:nlnn'''}
&&\sup_t\int_{B_1} n (|\ln n|+1) \zeta + 2 \int_{{Q_1}} |\nabla \sqrt{n}| \zeta + 2 \sup_t\int_{B_1} |\nabla \sqrt{\tilde{c}}|^2 \zeta\nonumber\\ \nonumber
&& \quad + \frac89 \int_{{Q_1}} |\nabla^2 \sqrt{\tilde{c}}|^2 \zeta+\sup_t\Lambda_0\int_{B_1} |u|^2 \zeta + \Lambda_0\int_{{Q_1}} |\nabla u|^2 \zeta\\ \nonumber
&\leq& C\Lambda_0\Lambda_1\varepsilon_2^{\frac13}.
\een
which and the assumption of $p$ imply the regularity due to Theorem \ref{thm:regularity-fixed}. Specially, we choose
\beno
C\Lambda_0\Lambda_1\varepsilon_2^{\frac13} \leq \frac{\varepsilon_1}{(\Lambda_0\Lambda_1)^{4+4\alpha_0}},
\eeno
which means
\beno
\varepsilon_2 \leq \frac{\varepsilon_1^3}{C(\Lambda_0\Lambda_1)^{15+12\alpha_0}}.
\eeno

For the second assumption (\ref{eq:condition-103}), if
\beno
(C')^\frac54 \varepsilon_2' \leq \frac{\varepsilon_1^\frac{15}4}{C(\Lambda_0\Lambda_1)^{\frac{75}4+15\alpha_0}},
\eeno
we have
\beno
\int_{Q_1} |n \ln n|^\frac32 \leq C' \varepsilon_2' + C' (\varepsilon_2')^\frac45 \leq \frac{\varepsilon_1^3}{C(\Lambda_0\Lambda_1)^{15+12\alpha_0}}.
\eeno

The proof is complete.


{\bf Proof of Corollary \ref{thm:3}.}
We will use a parabolic version of the Vitali covering lemma:
Let $\{J={Q_{z_\alpha,r_\alpha}}\}_\alpha$ be any collection of parabolic cylinders contained in a bounded subset of $\mathbb{R}^{4}$, and noting $J=J_x\times J_t$, there exist disjoint $Q_{z_j,r_j}\in J, j\in N$, such that any cylinder in $J$ is contained in $Q_{z_j,5r_j}$ for some $j$.

Letting
\beno
Q^{*}((x,t),r)=B(x,r)\times(t-\frac78r^2,t+\frac18r^2),
\eeno
it is a translation in time of $Q((x,t),r)$. Besides, $Q(z,\frac r2)\in Q^{*}(z,r)$.
Let $\mathcal{S}_R =\mathcal{S} \cap R$ for any compact set $R \subset Q_\frac12$. Fix any $\delta > 0$. 
Assume that for any $z_j = (x_j,t_j) \in \mathcal{S}_R$, by Theorem \ref{thm:lin}, there exists $0 < r_{z_j}=r_j < \frac\delta{10}$ such that 
\beno
\left( \int_{Q_{r_j}(z_0)} n^{\frac53} + |\nabla \sqrt{c+1}|^{\frac{10}{3}} + |u|^{\frac{10}{3}}+|p|^\frac53\right)\geq \frac12 r_{j}^{\frac53} \varepsilon_2' = \frac12 r_{j}^{\frac53} \frac{\varepsilon_1^\frac{15}4}{C(\Lambda_0^j\Lambda_1^j)^{\frac{75}4+15\alpha_0}}
\eeno
due to scaling. Here, $\Lambda_0^j = 108(\|c\|_{L^\infty((-r_j^2,0)\times \mathbb{R}^3)} + 1) \leq \Lambda_0$, and $\Lambda_1^j = (r_j \|\nabla \phi\|_{L^\infty((-r_j^2,0)\times \mathbb{R}^3)} + 1) \leq \Lambda_1$.
Thus,
\beno
\left( \int_{Q_{r_j}(z_0)} n^{\frac53} + |\nabla \sqrt{c+1}|^{\frac{10}{3}} + |u|^{\frac{10}{3}}+|p|^\frac53\right)\geq \frac12 r_{j}^{\frac53} \varepsilon_4,
\eeno
where $\varepsilon_4 = \frac{\varepsilon_1^\frac{15}4}{C(\Lambda_0\Lambda_1)^{\frac{75}4+15\alpha_0}}$ is independent of $j$.

Then,
\beno
\mathcal{S}_R \subset \bigcup_{j\in N} Q^{\ast}(z_j,2r_{z_j}).
\eeno
Let $r_j = r_{z_j}$ and $\{Q_{r_j}(z_j)\}_{j\in \mathbb{N}}$ be the countable disjoint subcover guaranteed by the Vitali covering lemma, then
\beno
\mathcal{S}_R \subset \bigcup_{j \in \mathbb{N}} Q^{\ast}(z_j, {10r_j})\quad 10r_j<\delta.
\eeno
Note that $Q({z_j, \frac {r_j}2})\in Q^{*}({z_j,{r_j}})$ are disjoint, then
\beno
\sum 10r_j^{\frac53} &\leq& \sum_j \frac{20}{\varepsilon_4} \left(\int_{Q_{r_j}(z_0)} n^{\frac53} + |\nabla \sqrt{c+1}|^{\frac{10}{3}} + |u|^{\frac{10}{3}}+|p|^\frac53\right).
\eeno
Since the bounded-ness of the right part, we have
\beno
\sum 10r_j^{\frac53} < + \infty,
\eeno
which means that $S_R$ has Labesgue measure $0$. Since the finite covering theorem, we know that for any open neighborhood $J = J_x \times J_t \subset {Q_1}$ of $S_R$ satisfies $Q_{r_z}(z) \subset J$,
\beno
\sum 5r_j^{\frac53} \leq \sum_j \frac{C}{\varepsilon_4} \left(\int_{Q_{r_j}(z_0)} n^{\frac53} + |\nabla \sqrt{c+1}|^{\frac{10}{3}} + |u|^{\frac{10}{3}}+|p|^\frac53\right).
\eeno
By the arbitrarily of $J$, we can choose $J$ with arbitrarily small Lebesgue measure, therefore, the right side is arbitrarily small. Since $\delta > 0$ is arbitrary, we have
\beno
{\mathcal{P}}^{\frac53}(\mathcal{S}_R) = 0.
\eeno
By the arbitrarily of $R$, we have
\beno
{\mathcal{P}}^{\frac53}(\mathcal{S}) = 0.
\eeno

\endproof

\section{Proof of Theorem {\ref{thm:2}}}

Next we prove Theorem \ref{thm:2} with the help of Theorem \ref{thm:regularity-fixed}.

{\it Proof:}
It's sufficient to prove the smallness of $p$ at a fixed ball, since other terms are the same as those in Theorem \ref{thm:regularity-fixed}.

{\bf Step I. The pressure estimate.} Recall the equation of $u$ in (\ref{eq:KS}),
and taking the divergence yields that
\beno
-\Delta p = \partial_i \partial_j (u_i u_j) + \nabla \cdot (n \nabla \phi).
\eeno
Let $\eta(x) \geq 0$ be supported in $B_\rho$ with $\eta = 1$ in $B_\frac\rho2$, and
\beno
p_1(x,t) = \int_{\mathbb{R}^3} \frac1{4\pi|x-y|} \left(\partial_i \partial_j ((u_i-(u_i)_\rho) (u_j-(u_j)_\rho) \eta) +\nabla (n\nabla\phi \eta)\right)(y,t) dy.
\eeno
Moreover, let
\beno
p_2(x,t) =p(x,t) - p_1(x,t)
\eeno
which implies that
\beno
\Delta p_2 = 0 ~~ {\rm in} ~~B_\frac\rho2.
\eeno
Let $0 < 2r < \rho \leq 1$, by the mean value property of harmonic functions and Lemma \ref{mean value property}, we have
\ben\label{ine:p2} \nonumber
\int_{{B_r}}|p_2|^{\frac32}dx &\leq& C\left(\frac r\rho\right)^{3} \int_{B_{\frac34 \rho}}|p_2|^{\frac32}dx \\
&\leq& C \left(\frac r\rho\right)^{3} \int_{B_{\rho}}|p|^{\frac32}dx+C \left(\frac r\rho\right)^{3}\int_{B_{\rho}}|p_1|^{\frac32}dx.
\een
and Calderon-Zygmund estimates yields that
\ben\label{ine:p1}
\int_{B_{\rho}}|p_1|^{\frac32}dx\leq C \int_{B_{\rho}}|u-(u)_\rho|^3+C \rho^{\frac 34}\left(\int_{B_{\rho}}|n\nabla \phi|^\frac65dx\right)^{\frac54}.
\een
Combining (\ref{ine:p1}) and (\ref{ine:p2}) and noting that $r<\rho$, we have
\beno
r^{-2} \int_{Q_r} |p|^\frac32 &\leq& r^{-2} \int_{Q_r} |p_1|^\frac32 + r^{-2} \int_{Q_r} |p_2|^\frac32 \\
&\leq& \left(1+\left(\frac r\rho\right)^{3}\right)r^{-2} \int_{Q_\rho} |p_1|^\frac32+ C r^{-2} \left(\frac r\rho\right)^{3} \int_{Q_{\rho}}|p|^{\frac32}dxdt \\
&\leq& C r^{-2} \int_{Q_{\rho}}|u-(u)_\rho|^3 + C r^{-2} \int_{I_\rho}\rho^{\frac 34}\left(\int_{B_{\rho}}|n\nabla \phi|^\frac65dx\right)^{\frac54}dt \\
&&+ C r^{-2} \left(\frac r\rho\right)^{3} \int_{Q_{\rho}}|p|^{\frac32}dxdt
\eeno
Using H\"{o}lder inequality, we have
\ben\label{pressure}
r^{-2} \int_{Q_r} |p|^\frac32
&\leq& C \left(\frac\rho r\right)^2 \rho^{-2} \int_{Q_{\rho}}|u-(u)_\rho|^3 + C \left(\frac\rho r\right)^2 \rho^\frac32 \left(\rho^{-\frac53} \int_{Q_\rho} |n|^\frac53dxdt\right)^\frac9{10} \|\nabla\phi\|_{L^\infty_{t,x}}^\frac32\nonumber\\
&&+ C \left(\frac r\rho\right) \rho^{-2} \int_{Q_{\rho}}|p|^{\frac32}dx.
\een

By Gagliardo-Nirenberg and Young inequality, noting that
\beno
r^{-1} \|n\|_{L^\frac53(Q_r)} \leq C r^{-1} \|\sqrt{n}\|_{L^\infty L^2 (Q_r)}^2 + C r^{-1} \|\nabla \sqrt{n}\|_{L^2(Q_r)}^2 \leq C (A_n(r) + E_n(r)),
\eeno
It follows from (\ref{pressure}) that
\ben\label{pressure2}
D(r) \leq C \left(\frac\rho r\right)^2 \tilde{C}_u(\rho) + C \left(\frac\rho r\right)^2 (A_n(\rho) + E_n(\rho))^\frac32 + C \left(\frac r\rho\right) D(\rho),
\een
where $\rho\leq 1$ is used.
Note that
\beno
\int_{B_{\rho}} |u-(u)_\rho|^3dxdt\leq \|u-(u)_\rho\|^\frac32_{L^2(B_\rho)} \|u-(u)_\rho\|^\frac32_{L^6(B_\rho)}\leq C\|u-(u)_\rho\|^\frac32_{L^2(B_\rho)} \|\nabla u\|^\frac32_{L^2(B_\rho)},
\eeno
and integrating in time we get
\beno
\rho^{-2}\int_{Q_{\rho}} |u-(u)_\rho|^3dxdt
&\leq&   C (A_u(\rho) + E_u(\rho))^\frac32.
\eeno
Then by (\ref{pressure2}) we have
\beno
D(r) \leq C \left(\frac r\rho\right) D(\rho) + C \left(\frac\rho r\right)^2 (A_n(\rho) + E_n(\rho))^\frac32+ C \left(\frac\rho r\right)^2 (A_u(\rho) + E_u(\rho))^\frac32.
\eeno

Let $G(\rho)=A_u(\rho) + E_u(\rho)+A_n(\rho) + E_n(\rho)$, we have
\beno
D(r) \leq C \left(\frac r\rho\right) D(\rho) + C \left(\frac\rho r\right)^2 G(\rho)^\frac32.
\eeno
Moreover, let $r = \theta_0 \rho$ with $\theta_0 \in (0,\frac14)$ satisfying $C\theta_0\leq \frac12$
and we arrive at
\beno
D(\theta_0 \rho) \leq \frac12 D(\rho) + C \theta_0^{-2} G(\rho)^\frac32.
\eeno
Write $\varepsilon'=\frac{\varepsilon_3}{(\Lambda_0\Lambda_1)^{4+4\alpha_0}}$.
Due to (\ref{eq:condition-all}), there exists $\rho_0>0$ such that $G(\rho_0) \leq 2\varepsilon'$. Then 
\beno
D(\theta_0 \rho_0) \leq \frac12D(\rho_0) + 4C \theta_0^{-2} \varepsilon'^{\frac32},
\eeno
and
\beno
&&D(\theta^2_0 \rho_0) \leq \frac12 D(\theta_0\rho_0) +C \theta_0^{-2} G(\theta_0\rho_0)^\frac32\leq \frac12 D(\theta_0\rho_0) +4 C \theta_0^{-2} \varepsilon'^{\frac32}\\
&&D(\theta^3_0 \rho) \leq \frac12D(\theta^2_0 \rho) +4 C \theta_0^{-2} \varepsilon'^{\frac32}\\
&& \quad \quad \quad \quad \quad \quad\cdots
\eeno
Consequently,  repeating the above progress, one can get
\beno
D(\theta^k_0\rho_0)\leq \frac1{2^{k}}D(\rho_0)+\sum_{k\geq 1}\frac1{2^{k-1}} (4C \theta_0^{-2} \varepsilon'^{\frac32}).
\eeno
which implies
\beno
D(\theta_0^k \rho_0) \leq \frac1{2^{k}}  D(\rho_0) + 8C \theta_0^{-2}\varepsilon'^{\frac32}.
\eeno
Recall that $\varepsilon_0=\frac{\varepsilon_1}{(\Lambda_0\Lambda_1)^{4+4\alpha_0}}$ in Theorem \ref{thm:regularity-fixed}. Nothing that $D(\rho_0) <\infty$, there exist a small constant $\varepsilon'$ satisfies  $\varepsilon'<\frac1{128}\varepsilon_0$ and $8C \theta_0^{-2}{\varepsilon'}^\frac32  \leq \frac{\theta_0^2}{32} \varepsilon_0$ and the constant $k = K_0(\theta_0)$ which is large enough and satisfies
\beno
2^{-K_0} D(\rho_0) \leq \frac{\theta_0^2}{32} \varepsilon_0,
\eeno
such that
\beno
D(\theta_0^{k} \rho_0) \leq \frac{\theta_0^2}{16}\varepsilon_0,
\eeno
holds for all $k \geq K_0$. Then for all $0<r\leq \theta_0^{K_0}\rho_0$, there exist constant $K_1 > K_0$ such that $\theta_0^{K_1+1} \rho_0 < r < \theta_0^{K_1} \rho_0$. Then
\beno
D(r) = r^{-2} \int_{Q_r} |p|^\frac32 \leq \theta_0^{-2K_1-2} \rho_0^{-2} \int_{Q_{\theta_0^{K_1} \rho_0}} |p|^\frac32 \leq \theta_0^{-2} \frac{\theta_0^2}{16}\varepsilon_0 \leq \frac1{16}\varepsilon_0.
\eeno

{\bf Step II. Arguments due to scaling.}
Recalling the assumption of (\ref{eq:condition-fixed}), we arrive at
\ben\label{eq:estimate-r}
&&\sup_t r^{-1}\int_{B_{r}} n  + |n \ln n| +|\nabla \sqrt{\tilde{c}}|^2 + |u|^2 \nonumber\\
&&+ r^{-1}\int_{Q_{r}} |\nabla \sqrt{n}|^2 + |\nabla u|^2 + |\nabla^2 \sqrt{\tilde{c}}|^2 + r^{-2}\int_{Q_{r}}|p|^\frac32 \leq \frac18\varepsilon_0,
\een
for any $0<r\leq r_0.$ By (\ref{eq:scaling}), let
\beno\label{eq:scaling'}
&&n_r(x,t)=r^2 n(r x, r^2t);~~ c_r(x,t)=c(r x, r^2t);\nonumber\\
&&~~ u_r(x,t)=r u(r x, r^2t);~~p_r(x,t)=r^2 p(r x, r^2t)
\eeno
then $(n_r,c_r,u_r,p_r)$ is a solution of (\ref{eq:KS}) in ${Q_1}^t$. Moreover,
\ben\label{eq:estimate-r2}
&&\sup_t \int_{B_{1}} n_r  +|\nabla \sqrt{\tilde{c_r}}|^2 + |u_r|^2 \nonumber\\
&&+ \int_{Q_{1}} |\nabla n_r^\frac12|^2 + |\nabla u_r|^2 + |\nabla^2 \sqrt{\tilde{c_r}}|^2 +|p_r|^\frac32 \leq \frac18\varepsilon_0.
\een
The remaining part is to estimate the term of $\sup_t \int_{B_{1}} |n_r\ln n_r|$, and
\beno
&&\sup_{-1<t<0} \int_{B_{1}} |n_r\ln n_r|\\
&\leq & \sup_{-r^2<t<0} r^{-1}\int_{B_{r}} |n\ln (r^2n)|dx\\
&\leq & \sup_{-r^2<t<0} r^{-1}\int_{B_{r}\cap\{n<r^{-\frac32}\}} |n\ln (r^2n)|dx+\sup_{-r^2<t<0} r^{-1}\int_{B_{r}\cap\{r^{-\frac32}\leq n\leq r^{-2}\}} |n\ln (r^2n)|dx\\
&&+\sup_{-r^2<t<0} r^{-1}\int_{B_{r}\cap\{n\geq r^{-2}\}} |n\ln (r^2n)|dx\doteq M_1+\cdots+M_3
\eeno
Then by (\ref{eq:estimate-r})
\beno
M_1\leq \sup_{-r^2<t<0} r^{-1}\int_{B_{r}\cap\{n<r^{-\frac32}\}} |n\ln n|+2n|\ln r|dx\leq \frac18\varepsilon_0+2r^{\frac12}|B_1||\ln r|\leq \frac14\varepsilon_0
\eeno
where we choose $r<r_1$ for a small $r_1 < r_0$.
Besides,
\beno
M_2&\leq&\sup_{-r^2<t<0} r^{-1}\int_{B_{r}\cap\{r^{-\frac32}\leq n\leq r^{-2}\}} |n\ln (r^2n)|dx\\
&\leq & \sup_{-r^2<t<0} r^{-1}\int_{B_{r}\cap\{r^{-\frac32}\leq n\leq r^{-2}\}} |n\ln (r^{-\frac12})|dx\\
&\leq & \sup_{-r^2<t<0} r^{-1}\int_{B_{r}} |n\ln n|dxdt\leq \frac18\varepsilon_0
\eeno
and
\beno
M_3&\leq&\sup_{-r^2<t<0} r^{-1}\int_{B_{r}\cap\{n\geq r^{-2}\}} |n\ln (r^2n)|dx\\
&\leq & \sup_{-r^2<t<0} r^{-1}\int_{B_{r}\cap\{n\geq r^{-2}\}} n\ln n dx\leq \frac18\varepsilon_0
\eeno
where we used $r<_1<1.$
Hence there holds $\sup_t \int_{B_{1}} |n_r\ln n_r|\leq \frac12 \varepsilon_0$. Combining this with (\ref{eq:estimate-r}) we get
\ben\label{eq:estimate-r2}\nonumber
&&\sup_t \int_{B_{1}} n_r +|n_r\ln n_r| +|\nabla \sqrt{\tilde{c_r}}|^2 + |u_r|^2 \\
&&+ \int_{Q_{1}} |\nabla n_r^\frac12|^2 + |\nabla u_r|^2 + |\nabla^2 \sqrt{\tilde{c_r}}|^2 +|p_r|^\frac32 \leq \varepsilon_0.
\een

Using Theorem (\ref{thm:regularity-fixed}), it follows that $(n_r,u_r,\nabla c_r)$ is regular at the point $0$. The proof is complete.

\noindent {\bf Acknowledgments.}
The author would like to thank Professors Sining Zheng and Zhaoyin Xiang for some helpful communications. W. Wang was supported by NSFC under grant 12071054,  National Support Program for Young Top-Notch Talents and by Dalian High-level Talent Innovation Project (Grant 2020RD09).


\begin{thebibliography}{WWW}




\bibitem{ABC} D. Albritton, E. Bru\'{e}, M. Colombo, {\it Non-uniqueness of Leray solutions of the forced Navier-Stokes equations.} Ann. of Math. (2) 196 (2022), no. 1, 415-455.


\bibitem{Bl2020} T. Black, {\it The Stokes limit in a three-dimensional chemotaxis-Navier-Stokes system.} J. Math. Fluid Mech. 22 (2020), no. 1, Paper No. 1, 35 pp.

\bibitem{BW2022} T. Black, M. Winkler, {\it Global weak solutions and absorbing sets in a chemotaxis-Navier-Stokes system with prescribed signal concentration on the boundary.} Math. Models Methods Appl. Sci. 32 (2022), no. 1, 137-173.

\bibitem{CKN} L. Caffarelli, R. Kohn and L. Nirenberg, {\it Partial regularity of suitable weak solutions of the Navier-Stokes equations}, Comm. Pure Appl. Math., 35(1982), 771-831.
\bibitem{CL2016} X. Cao, J. Lankeit, {\it Global classical small-data solutions for a three-dimensional chemotaxis Navier-Stokes system involving matrix-valued sensitivities,} Calc. Var. Partial Differential Equations 55 (2016), no. 4, Art. 107, 39 pp.



\bibitem{CKL2013} M. Chae, K. Kang, J. Lee, {\it On Existence of the smooth solutions to  the
coupled chemotaxis-fluid equations,} Discrete Cont. Dyn. Syst. A (2013), 33:2271-2297.

\bibitem{CKL2014} M. Chae, K. Kang, J. Lee, {\it Global existence and temporal decay in Keller-Segel models coupled to fluid equations,} Comm. Partial Differential Equations 39 (2014), no. 7, 1205-1235.


\bibitem{CKK2014} Y. Chung, K. Kang, J. Kim, {\it Global existence of weak solutions for a Keller-segel-fluid model with nonlinear diffusion,} J. Korean Math. Soc. 51 (2014) 635-654.




\bibitem{DL2020} M. Dai, H. Liu, {\it Low modes regularity criterion for a chemotaxis-Navier-Stokes system.} Commun. Pure Appl. Anal. 19 (2020), no. 5, 2713-2735.

\bibitem{DL2022} M. Ding, J. Lankeit, {\it Generalized solutions to a chemotaxis-Navier-Stokes system with arbitrary superlinear degradation.} SIAM J. Math. Anal. 54 (2022), no. 1, 1022-1052.

\bibitem{DCCGK2004} C. Dombrowski, L. Cisneros, S. Chatkaew, R.E. Goldstein, J.O. Kessler, {\it Self-concentration and large-scale coherence in bacterial dynamics,}
Phys. Rev. Lett. (2004), 93:098103-1-4.

\bibitem{DD} H. Dong and D. Du, {\it Partial regularity of solutions to the four-dimensional Navier-Stokes equations at the first blow-up time},
Comm. Math. Phys., 273(2007), 785-801.

\bibitem{DLX2017} R. Duan, X. Li, Z. Xiang, {\it Global existence and large time behavior for a two-dimensional chemotaxis-Navier-Stokes system.} J. Differential Equations 263 (2017), no. 10, 6284-6316.

\bibitem{DLM2010} R. Duan, A. Lorz, P. Markowich, {\it Global solutions to the coupled chemotaxis-fluid equations,} Comm. Partial Differential Equations 35 (2010) 1635-1673.



\bibitem{FLM2010} M. Francesco, A. Lorz and P. Markowich,{\it Chemotaxis-fluid coupled model for swimming bacteria with nonlinear diffusion: global existence and asymptotic behavior,} Discrete Continuous Dynam. Systems-A, 28 (2010), 1437-1453.

\bibitem{GKT} S. Gustafson, K. Kang and T.-P. Tsai, {\it Interior regularity criteria for suitable weak solutions of the Navier-Stokes equations,} Comm. Math. Phys., 273(2007), 161-176.






\bibitem{HZ2017} H. He, Q. Zhang, {\it Global existence of weak solutions for the 3D chemotaxis-Navier-stokes equations,} Nonlinear Analysis: Real World Applications. 35(2017) 336-349.

\bibitem{He2020} F. Heihoff, {\it Global mass-preserving solutions for a two-dimensional chemotaxis system with rotational flux components coupled with a full Navier-Stokes equation.} Discrete Contin. Dyn. Syst. Ser. B 25 (2020), no. 12, 4703-4719.

\bibitem{KLW2022}   K. Kang, J. Lee, M. Winkler, {\it Global weak solutions to a chemotaxis-Navier-Stokes system in $\mathbb{R}^3$.} Discrete Contin. Dyn. Syst. 42 (2022), no. 11, 5201-5222.

\bibitem{KP2002}
Katz, N. H.; Pavlovi\'{c}, N. {\it A cheap Caffarelli-Kohn-Nirenberg inequality for the Navier-Stokes equation with hyper-dissipation.} Geom. Funct. Anal. 12 (2002), no. 2, 355-379.

\bibitem{KM2019} S. Kurima, M. Mizukami, {\it Global weak solutions to a 3-dimensional degenerate and singular chemotaxis-Navier-Stokes system with logistic source.} Nonlinear Anal. Real World Appl. 46 (2019), 98-115.

\bibitem{LS} O. Ladyzhenskaya and G. Seregin, {\it On partial regularity of suitable weak solutions of the three-dimensional Navier-Stokes equations}, J. Math.Fluid.Mech., 1(1999), 357-387.

\bibitem{LR} P.  Lemari\'{e}-Rieusset. Recent developments in the Navier-Stokes problem, volume 431 of
Chapman \& Hall/CRC Research Notes in Mathematics. Chapman \& Hall/CRC, Boca Raton,
FL, 2002.

\bibitem{Leray} J. Leray, {\it Sur le mouvement d'un liquide visqueux emplissant l'espace.} Acta Math., 63(1):193-248, 1934.



\bibitem{LL2016} Y. Li, Y. Li, {\it Global boundedness of solutions for the chemotaxis-Navier-Stokes system in $R^2$.} J. Differential Equations 261 (2016), no. 11, 6570-6613

\bibitem{Lin} F. Lin, {\it A new proof of the Caffarelli-Kohn-Nirenberg theorem}, Comm. Pure Appl. Math., 51(3)(1998), 241-257.

\bibitem{LL2011} J. Liu  and  A. Lorz,   {\it A coupled chemotaxis-fluid model: Global existence,} Ann. Inst. H. Poincar\'{e} Anal. Nonlineaire, 28(2011), 643-652.



\bibitem{AL2010} A. Lorz, {\it Coupled chemotaxis fluid model,} Math. Models Methods Appl. Sci. 20 (2010) 987-004.




\bibitem{SV1} V. Scheffer,  {\it Partial regularity of solutions to the Navier-Stokes equations},  Pacific J. Math., 66(1976), 535-562.

\bibitem{SV2} V. Scheffer, {\it Hausdorff measure and the Navier-Stokes equations}, Commun. Math. Phy., 55(1977), 97-112.


\bibitem{SV4} V. Scheffer, {\it The Navier-Stokes equations on a bounded domain}, Commun. Math. Phy., 71(1980), 1-42.




\bibitem{Se} G.  Seregin, {\it Estimate of suitable solutions to the Navier-Stokes equations in critical Morrey spaces},
Journal of Mathematical Sciences, 143(2007),  2961-2968.







































\bibitem{TW2012} Y. Tao and M. Winkler, {\it Global existence and boundedness in a Keller-Segel-Stokes model with arbitrary porous medium diffusion,} Discrete Continuous Dynam. Systems-A, 32 (2012), 1901-1914.

\bibitem{TW2013} Y. Tao and, M. Winkler, {\it Locally bounded global solutions in a three-dimensional chemotaxis-Stokes system with nonlinear diffusion,} Ann. Inst. H. $Poincar\acute{e}$ Anal. Non $Lin\acute{e}aire$ 30 (2013) 157-178.



\bibitem{TX} G. Tian and Z. Xin, {\it Gradient estimation on Navier-Stokes equations}, Comm. Anal. Geom., 7(1999), 221-257.


\bibitem{Tsai-2018} Tsai, Tai-Peng, {\it Lectures on Navier-Stokes equations.} Graduate Studies in Mathematics, 192. American Mathematical Society, Providence, RI, 2018. xii+224 pp. ISBN: 978-1-4704-3096-2


\bibitem{TCDWKG2005} I. Tuval, L. Cisneros, C. Dombrowski, C. Wolgemuth, J. Kessler, R. Goldstein, {\it Bacterial swimming and oxygen transport near constant lines,} Proc. Natl. Acad. Sci. USA 102 (2005) 2277-282.








\bibitem{Va} A. Vasseur, {\it A new proof of partial regularity of solutions to Navier-Stokes
equations,} NoDEA Nonlinear Differential Equations Appl., 14 (2007), no. 5-6,
753-785.




\bibitem{WWX2018} Y. Wang, M. Winkler, Z. Xiang, {\it Global classical solutions in a two-dimensional chemotaxis-Navier-Stokes system with subcritical sensitivity.} Ann. Sc. Norm. Super. Pisa Cl. Sci. (5) 18 (2018), no. 2, 421-466.

\bibitem{WWX2018-2} Y. Wang, M. Winkler, Z. Xiang, {\it The small-convection limit in a two-dimensional chemotaxis-Navier-Stokes system.} Math. Z. 289 (2018), no. 1-2, 71108.

\bibitem{WZZ2021}  W. Wang, M. Zhang, S. Zheng, {\it To what extent is cross-diffusion controllable in a two-dimensional chemotaxis-(Navier-)Stokes system modeling coral fertilization.} Calc. Var. Partial Differential Equations 60 (2021), no. 4, Paper No. 143, 29 pp.


\bibitem{Winkler2012} M. Winkler, {\it Global large-data solutions in a chemotaxis-(Navier-)Stokes system modeling cellular swimming in fluid drops,} Comm. Partial Differential Equations 37 (2012) 319-352.

\bibitem{Winkler2016}  M. Winkler, {\it Global weak solutions in a three-dimensional chemotaxis-Navier-stokes system,} Ann. Inst. H. Poincar\'{e} C Anal. Non Lin\'{e}aire, 33 (2016), no. 5, 1329-1352.

\bibitem{Winkler2017} M. Winkler, {\it How far do chemotaxis-driven forces influence regularity in the Navier-Stokes system?} Trans. Amer. Math. Soc. 369 (2017), no. 5, 3067-3125.


\bibitem{ZZ2014} Q. Zhang, X. Zheng, {\it Global well-posedness for the two-dimensional incompressible chemptaxis-Navier-tokes equations,} SIAM J. Math. Anal. 46 (2014) 3078-3105.




\bibitem{ZK2021}
J. Zheng, Y. Ke, {\it Global bounded weak solutions for a chemotaxis-Stokes system with nonlinear diffusion and rotation.} J. Differential Equations 289 (2021), 182-235.























































%
%





%
%



%



%


%







%










%

%
















%
%
%
%







 \end{thebibliography}
\end{document}